\documentclass{article}
\usepackage[utf8]{inputenc}
\usepackage{fullpage}
\usepackage{hyperref}
\usepackage{amsfonts}
\usepackage{verbatim}
\usepackage{epsfig,graphics,color}

\usepackage{lineno,hyperref}
\modulolinenumbers[1]

\usepackage{array}
\newcolumntype{C}[1]{>{\centering\arraybackslash}m{#1}}

\usepackage{graphics}
\usepackage{amssymb, latexsym, amsmath, epsfig}
\usepackage{graphicx}
 \usepackage{color}
 \usepackage{epstopdf}
 \usepackage{comment}
 \usepackage{soul}
 \usepackage[normalem]{ulem}
  \usepackage{float}

\usepackage{amssymb,amsmath}
\usepackage{amsfonts,amstext,amsthm,amssymb,comment}
\usepackage{epsfig,graphics,color}

\graphicspath{{figures/}}

\usepackage{MnSymbol}

\numberwithin{equation}{section}

\usepackage{algorithm}
\usepackage{algorithmicx, algpseudocode}

\newcommand{\rev}[1]{{\color{black}#1}}
\newcommand{\ale}[1]{{\color{magenta} #1}}

\definecolor{samcolor}{rgb}{1, 0, 0}

\def\d{{\, \rm d}}
\newcommand{\bfs}[1]{{\boldsymbol #1}}

\title{Particle-Continuum Multiscale Modeling of Sea Ice Floes\thanks{The research of S.N.S and N.C. was partially funded by the Office of Naval Research (ONR) Multidisciplinary University Initiative (MURI) award N00014-19-1-2421. Q.D. was supported as a Van Vleck Visiting Assistant Professor under this grant. Computational resources for this project have been partially supported by the National Computing Infrastructure (NCI) national facility under Grant zv32.}}
\author{Quanling Deng\thanks{Corresponding author. School of Computing, Australian National University, Canberra, ACT 2601, Australia (Quanling.Deng@anu.edu.au). }
\and Samuel N. Stechmann\thanks{Department of Mathematics and Department of Atmospheric and Oceanic Sciences, University of Wisconsin--Madison, Madison, WI 53706, USA (Stechmann@wisc.edu)}
\and Nan Chen\thanks{Department of Mathematics, University of Wisconsin--Madison, Madison, WI 53706, USA
  (Nan.Chen@wisc.edu)}
  }

\begin{document}

\maketitle

\begin{abstract}
Sea ice profoundly influences the polar environment and the global climate. Traditionally, sea ice has been modeled as a continuum under Eulerian coordinates to describe its large-scale features, using, for instance, viscous-plastic rheology. Recently, Lagrangian particle models, also known as the discrete element method (DEM) models, have been utilized for characterizing the motion of individual sea ice fragments (called floes) at scales of 10 km and smaller, especially in marginal ice zones. This paper develops a multiscale model that couples the particle and the continuum systems to facilitate an effective representation of the dynamical and statistical features of sea ice across different scales. The multiscale model exploits a Boltzmann-type system that links the particle movement with the continuum equations. For the small-scale dynamics, it describes the motion of each sea ice floe. Then, as the large-scale continuum component, it treats the statistical moments of mass density and linear and angular velocities. The evolution of these statistics affects the motion of individual floes, which in turn provides bulk feedback that adjusts the large-scale dynamics. Notably, the particle model characterizing the sea ice floes is localized and fully parallelized, in a framework that is sometimes called superparameterization, which significantly improves computation efficiency. Numerical examples demonstrate the effective performance of the multiscale model. Additionally, the study demonstrates that the multiscale model has a linear-order approximation to the truth model.

\end{abstract}

\section{Introduction}

\begin{comment}
\begin{color}{red}
Some mathematicians who have done work on super-parameterization
(not necessarily sea ice):
\begin{enumerate}
    \item Yoonsang Lee
    \item Ian Grooms
    \item Stamen Dolaptchiev?
    \item Rupert Klein (particle--continuum, but not super-parameterization?)
\end{enumerate}
\end{color}
\end{comment}

Sea ice profoundly influences the polar environment and the global climate. Over the past several decades, there has been a significant sea ice loss in extent, thickness, and mass, possibly related to global warming and regional climate change \cite{mori2019reconciled, ogawa2018evaluating}. Therefore, it is crucial to model and analyze sea ice dynamics, which also facilitates developing the next-generation global climate models (GCMs) \cite{golden2020modeling, flato2004sea}. There has been a rapid growth in research on modeling sea ice, including contributions from applied and computational mathematics (see, for example, some recent work \cite{golden2020modeling, blockley2020future, roach2019advances, giddy2021stirring, chen2021lagrangian, manucharyan2022heavy, chen2022efficient, covington2022bridging, li2022improving, shih2023robust}).

%\begin{figure}[h!]
%\centering
%\hspace{-0.8cm}
%	\includegraphics[width=7in]{scale.png}
%	\caption{Multiscale structure of sea ice. }
%	\label{fig:ms}
%\end{figure}

% image comes from https://eoimages.gsfc.nasa.gov/images/imagerecords/83000/83962/hudsonbay_amo_2014169_lrg.jpg
% https://www.nasa.gov/sites/default/files/thumbnails/image/banner_0.jpg

\rev{Sea ice exhibits multiscale features \cite{golden2020modeling, hopkins2004formation}. The dynamics of sea ice across different scales interact and bring significant challenges to modeling sea ice, as described next. 
The traditional sea ice models focused on large scales ($\geq100$km) and treated sea ice as materials or plastics \cite{coon1974modeling, hibler1979dynamic, hunke1997elastic}. These models are based on the continuum assumption and capture the essential features such as ice thickness, concentration, deformation, and transport \cite{kreyscher2000results, docquier2017relationships, tandon2018reassessing}. 
Hibler's work \cite{hibler1979dynamic} in 1979 assumed that sea ice behaves as a plastic material
and developed viscous-plastic rheology. Later, the result was extended to an elastic version \cite{hunke1997elastic}. Alternative rheological models have also been developed to capture more features, such as fractures, by introducing a structure tensor variable \cite{wilchinsky2006modelling, dansereau2016maxwell}. These continuum models were helpful on large scales and were adopted in earth system models (ESMs); see, for example, \cite{stocker2014climate}.
These sea ice rheologies are continuum models with no concept of floe particles or sizes.
}

On the other hand, on scales of approximately $1$ to $10$ kilometers, the continuum assumption starts to break down \cite{holt2017prospects, coon2007arctic, feltham2008sea}, although in some cases, \rev{the continuum models are still valuable \cite{blockley2020future,kwok2008variability, kwok2002seasonal, hutter2020feature}
(for example, \cite{blockley2020future} states that these continuum models remain useful for large-scale and low-resolution modelling of sea ice as they cannot be readily invalidated using observation-based metrics).}
On such scales, the discrete element method (DEM) has been recognized as a promising alternative to modeling sea ice. Specifically, the DEM models are utilized to describe the motion of each sea ice fragment (called a floe) \cite{hopkins2004discrete, hopkins2004formation, herman2016discrete, bateman2019simulating, damsgaard2018application, xu2012discrete, tuhkuri2018review, manucharyan2022subzero} that better characterize granular and discontinuous materials \cite{cundall1979discrete, ghaboussi1990three, pande1990numerical, mishra1992discrete}. One notable difference between the continuum models and the DEMs is that the latter focuses on the physical quantities of sea ice floes along their trajectories under Lagrangian coordinates. Therefore, the DEM models are suited to describe the ice floes in marginal sea ice zone (MIZ) or where \rev{floes are fractured into smaller floes \cite{sammonds1998fracture}}. \rev{MIZ is the transition area between the open ocean and the pack ice, characterized by the presence of fragmented or thin ice floes. In the MIZ, modeling ice drift, breakup, and formation is particularly challenging. Many physical processes in the MIZ occur at scales smaller than the model's grid resolution, which brings challenges in parameterizations as well as in characterizing the floe sizes (see floe size distributions, for example, in \cite{alberello2022three,roach2018emergent,meylan2021floe}). Understanding the interaction between ocean waves and sea ice in MIZs is crucial for accurately modeling and predicting sea ice dynamics, especially in the context of climate change and its impact on polar regions;  we refer to  \cite{squire1996moving,squire2022marginal,dumont2022marginal,bennetts2022marginal,herman2022granular} and the references therein for details.  The DEM particle models} are especially appealing for operational forecasting applications that require models to reproduce sea ice behavior on fine spatiotemporal scales.

Despite these merits, there exist challenges in developing and simulating DEM models. First, special attention is needed to couple the Lagrangian DEM model with the traditional ocean and atmosphere dynamics under the Eulerian coordinates. Next, as the DEM model focuses on characterizing the dynamics of each ice floe, the simulation requires extensive computational resources. In addition, the parameterization of the floe dynamics to the continuum quantities may become inaccurate when the floe displays multiscale features. Hence, the DEM model has yet to be widely incorporated into GCMs.

\rev{This paper develops a multiscale sea ice model that combines the large-scale continuum dynamics with the small-scale Lagrangian DEM \cite{delle2020particle, delgado2003continuum, petsev2017coupling}. In this multiscale modeling framework, a so-called super-parameterization \cite{grabowski1999crcp, randall2003breaking, majda2009simple, stechmann2014multiscale,  grooms2013efficient, majda2014new, lee2017stochastic, dolaptchiev2013multiscale, jansson2019regional} is applied to the Lagrangian floe models.
One of the new challenges for superparameterization for sea ice floes, in comparison to past superparameterization work, is that the fine-scale model for sea ice is not an Eulerian continuum model but a Lagrangian particle (floe) model. To aid in connecting the fine-scale Lagrangian particle (floe) model and a coarse-scale Eulerian continuum model, a Boltzmann-type equation is used here \cite{cercignani1988boltzmann, harris2004introduction, de1989incompressible, golse2004navier}. In this particle--continuum superparameterization, an independent DEM model is adopted in each sub-domain to capture the particle features at fine scales and is dynamically coupled with the coarse scales. Since the different DEMs are evolving independently (in between the fine--course coupling that occurs on coarse time scales), \rev{the computations are fully} parallel in a parallel computing paradigm and help facilitate computational efficiency.

The primary mathematical contribution of the study lies in the integration of two scales in Eulerian and Lagrangian coordinates through a superparameterization technique which presents a significant advancement in understanding and simulating sea ice dynamics.
Traditional DEM-based sea ice models typically focus solely on the small-scale dynamics of individual ice floes, neglecting the interactions and feedback with the larger-scale continuum dynamics as they are local dynamics in Lagrangian coordinates. In contrast, the introduced multiscale model bridges this gap by incorporating both fine-scale particle features and coarse-scale continuum dynamics through a Boltzmann-type equation.
The proposed multiscale model brings a new perspective to sea ice research, offering improved insights into the complex interactions between different scales
and
making it a promising tool for studying sea ice dynamics in large-scale scenarios. 
}

This paper is one of few in the growing area of multiscale models for sea ice. \rev{The present approach is set up for an idealized DEM particle model. It serves as a multiscale framework designed to facilitate the future use of more complex floe models \cite{hopkins2004discrete, hopkins2004formation, herman2016discrete, bateman2019simulating, damsgaard2018application, xu2012discrete, tuhkuri2018review, manucharyan2022subzero} and coupling with atmospheric multi-scale models \cite{grabowski1999crcp, randall2003breaking, jansson2019regional}.
A different approach is taken by \cite{davis2023super} that facilitates the mathematical analysis and utilizes graph-theoretic techniques \cite{berry2016variable,davis2023graph}.
Another work in progress is to use a hybrid approach where a continuum 
model is solved in the whole domain, and a DEM is nested into the continuum model in regions of interest (e.g. MIZ) 
where the continuum approximation begins to break down
(Carolin Mehlmann, personal communication, 2023).
This progressive research underscores the ongoing trend towards multiscale modeling of sea ice and floe dynamics.}

%instead of using Lagrangian trajectories with fine-scale locations $(x^l(t),y^l(t))$ and velocities $(u^l(t),v^l(t))$, only the fine-scale velocities are used in the method of \cite{davis2023super}.

%The method results in a set of continuum equations governing the statistical moments of the sea ice's concentration (or mass density), linear velocity, and angular velocity. A Boltzmann-type equation is obtained to describe each of the continuum models \cite{cercignani1988boltzmann, harris2004introduction, de1989incompressible, golse2004navier}. These continuum equations aim to capture large-scale features. On the other hand, an independent DEM model with a reduced computational cost is adopted in each sub-domain to capture the particle features at fine scales. The fine-scale DEM models create statistical inputs, such as the floe velocity field, for the continuum equations. In contrast, the large-scale solutions provide mass density (concentration) and large-scale averaged momentum to update the individual velocity and sizes of the fine-scale floes. Notably, with a suitable domain decomposition, the model allows a fully parallel computing paradigm for the DEM models and thus facilitates computational efficiency.

The rest of this paper is organized as follows. Section~\ref{sec:md} starts by introducing the simplified DEM for sea ice floes, followed by super-parameterizing the DEM particle model using the particle-continuum model. A Boltzmann-like description is utilized for the floe particles, which characterize the time evolution of the probability distribution function in terms of the physical quantities of floes. Integrating the Boltzmann equation over hyperplanes arrives at the lower-order statistical moments equations that approximate the particle model's large-scale features. Section \ref{sec:alg} presents the methods for solving the particle model and the continuum equations, followed by an algorithm that details the coupling of the two sets of models. Section~\ref{sec:num} collects numerical results demonstrating the performance of the proposed method. Several simulation scenarios are included. The linear convergence of the $L^2$-norm errors of the floe concentration is also shown in this section. Concluding remarks are presented in Section~\ref{sec:con}.

\section{The particle-continuum model} \label{sec:md}

In this section, \rev{we first present a simple DEM particle model for sea ice floes following recent work \cite{herman2016discrete,damsgaard2018application, chen2021lagrangian,chen2022superfloe}.
With the particle model setting, assuming the particle number goes to infinity, we derive the Boltzmann-type equation and their moments. 
In particular, we adopt a continuum Boltzmann model based on the distribution function with a phase space consisting of radius, position, angular position, velocity, and angular velocity of each sea ice floe. Then the equations of lower-order moments are derived, which capture the large-scale features of the particles and serve as the coarse-scale continuum model.
In the numerical experiment section \ref{sec:num}, we use this DEM as a benchmark model for comparison with the proposed multiscale model. 
}

\subsection{DEM particle model} \label{sec:dem}

For the simplicity of illustrating the main idea of developing the multiscale model, a simple DEM model describing each sea ice floe is presented here. In this DEM model, the geometry of the sea ice floes is assumed to be all cylinders with \rev{a fixed unit thickness. In reality, the floe thickness distribution usually follows a Gamma distribution (see, for example, \cite{bourke1987sea, toppaladoddi2015theory} for the Arctic region and \cite{toyota2011size} for the Antarctic region). }Denote by $r^l$ the radius of the $l$-th floe with index $l=1,2, \cdots, N$, where $N$ is the total number of floes in the system. The radius characterizes the floe size, which generally follows a power law distribution \cite{stern2018seasonal}. With the unit thickness, the mass of the $i$-th floe is $m^l=\rho_{ice}\pi (r^l)^2 $, where the constant $\rho_{ice}$ is the density of sea ice floes. The floe position is denoted by $\bfs{x}^l=(x^l,y^l)^\mathtt{T}$ and the angular location is $\theta^l$. Similarly, the floe velocity is $\bfs{v}^l=(u^l,v^l)^{\mathtt{T}}$ and the angular velocity is $\boldsymbol{\omega}^l = \omega^l\hat{\bfs{z}}$. Herein, $\hat{\bfs{z}}$ is the unit vector along the $z$-axis (perpendicular to the $(x,y)$ plane). Lastly, let $\bfs{u}_o$ be the ocean surface velocity, which is assumed to be known.
With these notations in hand, a simple DEM system describing the sea ice floe dynamics is as follows:
\begin{subequations}\label{eq:dem}
\begin{align}
  \frac{\d\boldsymbol{x}^l}{\d t} &= \boldsymbol{v}^l,\label{dem_x}\\
  \frac{\d\theta^l}{\d t}  &= \omega^l,\label{dem_Omega}\\
  m^l\frac{\d\boldsymbol{v}^l}{\d t} &=    \sum_{j=1}^{\rev{N}}(\boldsymbol{f}_{\boldsymbol{n}}^{lj} + \boldsymbol{f}_{\boldsymbol{t}}^{lj}) + \widetilde{\alpha}^l\left(\boldsymbol{u}^l_o -\boldsymbol{v}^l\right)\left| \boldsymbol{u}^l_o -\boldsymbol{v}^l\right|
  := \boldsymbol{F},\label{dem_v}\\
  I^l\frac{\d\omega^l}{\d t} &= \sum_{j=1}^{\rev{N}} ( r^l\boldsymbol{n}^{lj}\times \boldsymbol{f}_{\boldsymbol{t}}^{lj}) \cdot\hat{\boldsymbol{z}} + \widetilde{\beta}^l \left(\nabla\times \boldsymbol{u}^l_o/2-\omega^l\hat{\boldsymbol{z}}\right)\left|\nabla\times \boldsymbol{u}^l_o/2-\omega^l\hat{\boldsymbol{z}}\right|
  := F_\omega,\label{dem_omega}
  %
  %\frac{\d\boldsymbol{\hat{u}}_o}{\d t} &= \big(\boldsymbol{L}_{\boldsymbol{u}}\boldsymbol{\hat{u}}_o + \boldsymbol{F}_{\boldsymbol{u}}\big) + \boldsymbol\Sigma_{\boldsymbol{u}} \dot{\boldsymbol{W}}_{\boldsymbol{u}}(t), \quad
  %\boldsymbol{u}^l_o = \boldsymbol{G}(\boldsymbol{x}^l) \hat{\boldsymbol{u}}_o, \label{dem_u}
\end{align}
\end{subequations}
where $l=1, 2, \cdots, N$.
\rev{The model equation \eqref{eq:dem} describes the dynamics of each particle. 
They are equations of motion following Newton’s law; we refer to \cite{chen2022superfloe,chen2021lagrangian} for details.
The floe contact force, being nonzero only when two floes are in contact, consists of the normal and tangential components, $\bfs{f}_\bfs{n}^{lj}$ and $\bfs{f}_\bfs{t}^{lj}$, respectively, 
where }
\begin{equation}
\begin{gathered}
     \boldsymbol{f}_{\boldsymbol{n}}^{lj} = c^{lj}E^{lj}\delta^{lj}\boldsymbol{n}^{lj},\qquad
     \boldsymbol{f}_{\boldsymbol{t}}^{lj} = c^{lj}G^{lj}v_{\boldsymbol{t}}^{lj}\boldsymbol{t}^{lj},\qquad
     \delta_n^{lj} \equiv d^{lj} - (r^l+r^j),
    \quad\mbox{with}\quad d^{lj} = |\bfs{x}^l-\bfs{x}^j|,\\
     v_{\boldsymbol{t}}^{lj} = \big[(\boldsymbol{v}^j + \boldsymbol{\omega}^j \times \boldsymbol{r}^j) - (\boldsymbol{v}^l + \boldsymbol{\omega}^l \times \boldsymbol{r}^l)\big]\cdot\boldsymbol{t}^{lj},\qquad
     \widetilde{\beta}^l = d_o \rho_o\pi (r^l)^4,\qquad
     \widetilde{\alpha}^l = d_o \rho_o\pi (r^l)^2.
\end{gathered}
\end{equation}
\rev{The first equation is Hooke's linear elasticity law for the normal force. 
Therein, $E^{lj}$ is Young's bulk modulus and $c^{lj}$ is the chord length in the transverse direction of the cross-sectional area. 
The tangential force describes the resistance against slip between floes by limiting relative tangential movement \cite{cundall1979discrete}. 
$G^{lj}$ is Young's shear modulus.
The vectors $\bfs{n}$ and $\bfs{t}$ are the unit vectors along the normal and the tangential directions,  and $\bfs{r}^j$ is the radius multiplied by the associated normal vector, defined by pointing towards the center of the $l$-th floe. 
The Coulomb friction law plays an important role in limiting the tangential contact force relative to the magnitude of the normal contact force \cite{hopkins2004discrete}. Thus, the restriction $|\bfs{f}_\bfs{t}^{lj}| \leq \mu^{lj}|\bfs{f}_\bfs{n}^{lj}|$ is applied, where $\mu^{lj}$ is the coefficient of friction that characterizes the condition of the surfaces of the two floes in contact. Assume that $\mu^{lj}$ is constant for all the floes. 
For the rest of the parameters, 
$d_o$ is the ocean drag coefficient, $\rho_o$ is the ocean density,
and $I^l=m^l(r^l)^2$ is the moment of inertia.
}

The model \eqref{eq:dem} neglects the atmospheric forcing for simplicity. \rev{In general, the atmospheric wind that induces forcing often contributes to the sea ice dynamics on large scales and can be a main driver of the sea ice motion (see, for example, \cite{alberello2020drift}). In such a case, one may include the wind drag force by adding similar terms to the ocean drag forces in \eqref{eq:dem}. 
We also assume no ridging, thus there are no out-of-plane motions in this DEM model. }
%With this in mind, we omit the atmospheric wind

It is worthwhile to highlight a few properties of the model. First, nonlinearity exists in the model, including the quadratic terms in the momentum and angular momentum equations, representing the ocean drag effect, and the more complicated nonlinearity in the contact forces. Second, the system is high-dimensional as each floe gives one set of state variables containing different physical quantities. Thus, when the floe number $N$ increases, the dimension of the system increases significantly. Lastly, although the model has idealized geometric properties, the coupled system \rev{captures many features of sea ice floe dynamics, making it a suitable test model for developing a multiscale system. 
We demonstrate this through a variety of numerical simulations. }

\subsection{Continuum description and distribution function}

\rev{The goal of this section is to develop a particle--continuum multiscale model. 
In this subsection, we briefly outline the motivation for the multiscale model development, followed by the introduction of a Boltzmann-type description of particles. }

\subsubsection{Motivation for particle--continuum super-parameterization}

\begin{figure}[h!]
    \centering
    \includegraphics[width=0.9\textwidth]{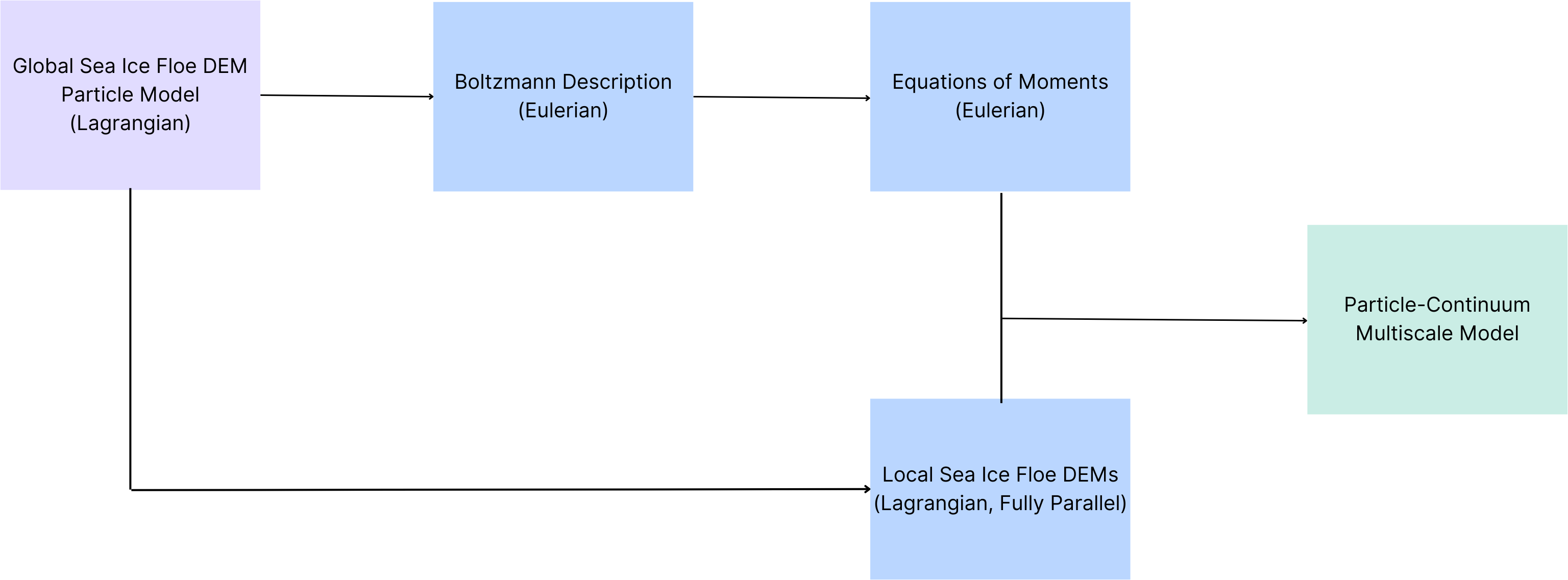}
    \caption{A diagram of particle--continuum multi-scale model.}
    \label{fig:modiag}
\end{figure}

One motivation for a hybrid particle--continuum \rev{model (see Figure \ref{fig:modiag} for a schematic of the model development)} of sea ice,
as opposed to the direct use of a particle model 
such as (\ref{eq:dem}) by itself, is for accurately representing
the fine scales of sea ice in multi-scale climate models.
Multi-scale models of the atmosphere are already in use in climate models
(and have been in use for 
roughly two decades
\cite{grabowski1999crcp,randall2003breaking,jansson2019regional}),
but sea ice is not typically represented via a multi-scale sea ice model. 
To have the sea ice and atmosphere and ocean represented 
on the same coarse and fine scales, 
it is natural to use a particle model for floes on
the fine scales. While particle models of sea ice were
less well-developed two decades ago during the advent of
atmospheric multiscale modeling (i.e., superparameterization), 
particle models of sea ice have undergone
substantial development in recent years
\cite{hopkins2004discrete, hopkins2004formation, herman2016discrete, bateman2019simulating, damsgaard2018application, xu2012discrete, tuhkuri2018review, manucharyan2022subzero},
and it is feasible to consider a multi-scale model of sea ice now
and to use a particle--continuum approach. 

%A large piece of sea ice floe can be approximated as a group of smaller floes in contact connecting with strong bonds. Within a short period, each floe in such a group has similar dynamical properties as the entire large piece of the floe. Yet, over an extended period, the large piece of the floe can fracture due to the stress tensor (essentially from the ocean and atmospheric drag and floe-contact forces). The breakdown can be approximately described by separating the smaller floes as a response to different drag and contact forces. The approximation becomes more accurate when more small floes are utilized. Hence, a sufficient number of ice floes is preferred. In regions where sea ice is individual pieces, such as the MIZ, it is ideal for representing each floe as an individual or a group of such particles to capture the associated dynamics behaviour. In the following, the entire sea ice field is treated as consisting of a large number of particles. The dynamics of individual particle floe follow the governing equations in \eqref{eq:dem}. However, in the presence of a large number of floes, it is generally computationally expensive to carry out numerical integration due to the necessity of using small time-step sizes to guarantee numerical stability and accuracy.

Another motivation for a hybrid particle--continuum model of sea ice
is the high computational cost of the direct use of a particle model 
such as (\ref{eq:dem}) by itself.
In the presence of a large number of floes, it is generally computationally expensive to carry out numerical integration due to the necessity of using small time-step sizes to guarantee numerical stability and accuracy.
One natural way to improve computational efficiency is to implement the simulation in a parallel environment directly. However, such a method may have the following potential issues. First, it remains challenging to optimally distribute the particles to different processors, which minimizes the message passing. Note that each floe particle may be in contact with neighboring floes. The state-of-the-art technique for parallelization in such a situation is to represent the contacting particles as graphs. However, graph partitioning is a challenging problem, especially in high-dimensional systems \cite{hendrickson2000graph,bulucc2016recent}. Secondly, the floes initially far from each other may come into contact later as the system evolves dynamically. This further requires graph partitioning to be dynamic, and the optimal partition would depend on the particle model. For instance, given a partition, certain floes near the partition boundary may bounce back and forth, requiring frequent floe data transfer between neighbouring processors. In addition, when passing one floe from one processor to another, certain unnecessary information transfers may occur. For example, in a one-time step during the dynamic process using parallel computing, one processor may receive new floes from its left processor and send floes to its right processor. These floes may have similar physical quantities: their floe velocities and radii or masses may be approximately equal. In this case, instead of transferring all the floe information among these three processors, it may be more efficient to apply periodic boundary conditions (floes leaving the right boundary to replace the ones entering from the left boundary) to each processor and only transfer the statistical moments (averages and variances) such as concentration and averaged momentum among the processors after several time steps (coarse time step).
\rev{In this paper, we propose a particle-continuum model multiscale framework to circumvent these challenges 
by fully parallelizing the fine-scale DEM particle models.}

Thus, an effective particle-continuum coupled model is needed to describe sea ice and its coupling with the ocean and atmosphere. The continuum model characterizes the large-scale feature of sea ice dynamics, and the DEM models the fine-scale parts. In addition, a suitable parallelization on the fine-scale model is highly desired to alleviate the computational burden.

%Based on the discussions above, we aim to super-parameterize the DEM model \eqref{eq:dem} to obtain a particle-continuum model. 

%To this end, we introduce the Knudsen number, denoted as $Kn$, which characterizes the scales of particles in the view of mechanics. The Knudsen number is a non-dimensional number defined as the ratio of the particle (originally for molecular) mean free path length to a representative physical length scale \cite{munafo2014multi}. For small values of $Kn$, one may use Boltzmann kinetic equation or Navier-Stokes-type equations to model the dynamics of particles under consideration. In contrast, for large values of $Kn$, one tends to adopt the Boltzmann equation for modelling the statistics of the particles. As sea ice exhibits multi-scale features, $Kn$ varies from very small values for large-scale floes to large values for small-scale floes in, for example, MIZs. When the Knudsen number is near or greater than one, the mean free path of a molecule is comparable to the length scale of the problem, and the continuum assumption of fluid mechanics is no longer valid. In such cases, statistical modelling methods using Boltzmann-type descriptions are preferred.

%\begin{figure}[h!]
%\centering
%	\includegraphics[width=6in]{knn.png}
	%\includegraphics[width=6in]{kn.jpg}
%	\caption{Dynamic regime classification according to the Knudsen number. }
%	\label{fig:kn}
%\end{figure}

%The foundation of the kinetic theory lies in the Boltzmann equation, which was initially introduced to describe the statistical behaviour of gas molecules 

\subsubsection{A Boltzmann description}

Based on the discussions above, we aim to super-parameterize the DEM model \eqref{eq:dem} to obtain a particle-continuum model. 
As a continuum description of the floe dynamics, we will use a kinetic theory formulation with the Boltzmann equation.
The Boltzmann equation was initially introduced to describe the statistical behavior of gas molecules and their velocity distribution
\cite{cercignani1988boltzmann,harris2004introduction}. In the following, the Boltzmann description is applied to characterize the statistical behavior of the sea ice DEM system with many floes. We will start with the Boltzmann equation that describes the probability distribution function (PDF) of the floes defined in a phase space consisting of radii, linear and angular positions, and linear and angular velocities. Yet, such a partial differential equation (PDE) contains multiple variables and is thus computationally expensive. To overcome such a computational challenge, the moment equations of sea ice floe dynamics are derived, which serve as the large-scale model when coupled with the local fine-scale DEM models.

Consider the phase space $S$ consisting of the set of all possible radii $r$, angular position $\theta$, angular velocities $\omega$, positions $\boldsymbol{x} = (x,y)$, and velocity $\boldsymbol{v}=(v_x, v_y)$ in two-dimensions (2D). Therefore, the entire phase space is seven-dimensional and an element in this space is
$(r, \omega, \theta, \boldsymbol{x}, \boldsymbol{v}) = (r, \omega, \theta, x, y, v_x, v_y) \in S$. By including the time $t$, the floe distribution function is defined as
$f(r, \omega, \theta, \boldsymbol{x}, \boldsymbol{v}, t) = f(r, \omega, \theta, x, y, v_x, v_y, t)$,
which gives the probability per unit phase-space volume  at an instance of time $t$.
We introduce the notation for the differential volume
\begin{equation}
    dE = dr \ d\omega \ d\theta \ d\boldsymbol{v} = dr \ d\omega \ d\theta \ dv_x\ dv_y
\end{equation}
and
\begin{equation}
    dE_a = \frac{dE}{da}, \quad a = r, \omega, \theta, \boldsymbol{v},
\end{equation}
which stands for the differential volumes over a subset of variables.
With these notations, the number density with respect to the 2D plane at time $t$ and the total number of particles are
\begin{equation} \label{eq:numdens}
    n = \int f \ dE, \quad \text{and} \quad N = \int n \ d\boldsymbol{x},
\end{equation}
respectively.
Other number density functions can be defined similarly.
For example, the number density with respect to the radius at time $t$ is
$
n_r = \int f \ dE_r.
$
The average of any function $g$ in a 2D plane is
\begin{equation} \label{eq:meang}
    \langle g \rangle = \frac{1}{n}\int g f \ dE,
\end{equation}
while the \rev{spatial} average is
\begin{equation}
    \overline{ \langle g \rangle } = \int \langle g \rangle \ d\boldsymbol{x}.
\end{equation}
Herein, the integration layers are complete.
In particular, the average radius is
\begin{equation}
    \langle r \rangle = \frac{1}{n} \int r f \ dE.
\end{equation}

Note that, different from the conventional phase space that consists of only position and velocity (or momentum), the phase space defined here is a generalized one. Below, we first describe the Boltzmann equation and then build a simplified model characterizing the associated moments to reduce the computational cost.

From the law of conservation, the general Boltzmann equation \rev{(see, for example, \cite{harris2004introduction})} reads
\begin{equation} \label{eq:gbe}
    \frac{\partial f}{\partial t}
    + \nabla_{\boldsymbol{x}} \cdot \Big(f \frac{d \boldsymbol{x}}{d t} \Big)
    + \nabla_{\boldsymbol{v}} \cdot \Big(f \frac{d \boldsymbol{v}}{d t} \Big)
    + \partial_\omega \Big(f \frac{d \omega}{d t} \Big)
    + \partial_\theta \Big(f \frac{d \theta}{d t} \Big)
    + \partial_r \Big(f \frac{d r}{d t} \Big)
    = C(f),
\end{equation}
where $C(f) = \Big( \frac{\partial f}{\partial t}\Big)_{\text{coll}}$
is collision contact term. 
\rev{The equation describes how the distribution function changes over time due to the advection of particles in phase space $S=(r, \omega, \theta, x, y, v_x, v_y)$ and the effects of particle collisions. It gives the probability of finding particles with a certain radius, position, velocity, angular position, and angular velocity at a given time.
The phase space for the standard Boltzmann equation consists of only particle positions and velocities \cite{harris2004introduction}). Our phase space is novel and has a larger dimension including variables describing particle size, angular position, and angular velocity. 
}
The above Boltzmann equation differs from the classical one for gas molecules that has an explicit expression for $C(f)$ under the two assumptions: (a) binary collisions are dominant; (b) molecular chaos. In the absence of collision, the Boltzmann equation is called the Vlasov equation. Applying to the sea ice floes, the only assumption is that the collisions happen instantly so that $C(f)$ satisfies the conservation law for number density, mass, momentum, and angular momentum. That is,
\begin{equation} \label{eq:gc0}
    \int g C(f) \ dE = 0,
\end{equation}
where $g = 1, m, m\boldsymbol{v}, I\omega$.

Omitting the superscripts in the DEM model \eqref{eq:dem} and plugging these equations into the Boltzmann equation \eqref{eq:gbe} arrive at the equation for floe dynamics
\begin{equation} \label{eq:dem4gbe}
    \frac{\partial f}{\partial t}
    + \nabla_{\boldsymbol{x}} \cdot ( \boldsymbol{v} f )
    + \nabla_{\boldsymbol{v}} \cdot \Big(f \frac{ \boldsymbol{F}}{m} \Big)
    + \partial_\omega \Big(f \frac{F_\omega}{I} \Big)
    + \partial_\theta (\omega f )
    = C(f),
\end{equation}
where the floe radius is \rev{assumed to be unchanged within the simulation, i.e., $\frac{dr}{dt} = 0.$ In reality, this term resembles partially the features of sea ice fracturing, freezing and melting, which are related to the sea ice surface temperature and salinity. No freezing nor melting occurring is assumed here for simplicity.}

%\subsection{Super-parameterization: equations of statistical moments}
\subsection{Statistical moments}
The equation \eqref{eq:dem4gbe} has a total of eight dimensions (including time); thus, it is very challenging to solve with standard discretization schemes. To this end, the equations of statistical moments are derived to reduce the space's dimension while capturing the physical quantities' main dynamical features on large scales. The macroscopic floe flow properties of interest can be derived by computing the moments of $f$ defined as integrals over all the phase spatial coordinates except $x$ and $y$, $dE$. The integrands are products of $f$ with different quantities $g$ carried by a floe; see \eqref{eq:meang}. As shown in equation \ref{eq:numdens}, the number density is computed as the moment with $g=1$. Below, we derive the moment equation from \eqref{eq:dem4gbe} starting with a generic physical quantity $g$.

We multiply equation \eqref{eq:dem4gbe} by an arbitrary function $g = g(r,\omega, \boldsymbol{v})$ (this can be more general but we restrict it here to cover the conserved quantities of interest) and take integration of both sides over $r, \omega, \theta,$ and $\boldsymbol{v}$. This leads to
\begin{equation}
    \int g \Big( \frac{\partial f}{\partial t}
    + \nabla_{\boldsymbol{x}} \cdot ( \boldsymbol{v} f )
    + \nabla_{\boldsymbol{v}} \cdot \Big(f \frac{ \boldsymbol{F}}{m} \Big)
    + \partial_\omega \Big(f \frac{F_\omega}{I} \Big)
    + \partial_\theta (\omega f ) \Big) \ dE
    = \int g C(f) \ dE.
\end{equation}
We now calculate it term by term to obtain the following.

\begin{itemize}
\item To calculate the first term, one applies the exchange of integration with time derivative and obtains
\begin{equation}
    \begin{aligned}
        \int g \frac{\partial f}{\partial t} \ dE
        = \frac{\partial }{\partial t} \int g f \ dE
        = \frac{\partial }{\partial t} \big(n \langle g \rangle \big).
    \end{aligned}
\end{equation}

\item To calculate the second term, we note that in the phase space, variables are independent. Thus, by exchange of integration with $\nabla_{\boldsymbol{x}}$, one arrives at
\begin{equation}
    \begin{aligned}
        \int g \nabla_{\boldsymbol{x}} \cdot (\boldsymbol{v} f) \ dE
        = \int \nabla_{\boldsymbol{x}} \cdot (g \boldsymbol{v} f ) \ dE
        = \nabla_{\boldsymbol{x}} \cdot \int \boldsymbol{v} g f  \ dE
        = \nabla_{\boldsymbol{x}} \cdot \big(n \langle g \boldsymbol{v} \rangle \big).
    \end{aligned}
\end{equation}

\item To calculate the third term, we first write the integration over $dE$ as double integration over $d\boldsymbol{v}$ and $dE_{\boldsymbol{v}}$.
By exchange of integration, integration by parts, and the assumption that $f$ is vanishing at the boundaries, one obtains
\begin{equation}
    \begin{aligned}
        \int g \nabla_{\boldsymbol{v}} \cdot \Big(f \frac{ \boldsymbol{F}}{m} \Big) \ dE
        & = \int \Big( \int g \nabla_{\boldsymbol{v}} \cdot \Big(f \frac{ \boldsymbol{F}}{m} \Big)  \ d\boldsymbol{v} \Big) \ dE_{\boldsymbol{v}} \\
        & = \int \Big( \int gf \frac{ \boldsymbol{F}}{m} \cdot \boldsymbol{n}_{\boldsymbol{v}} \ dS_{\boldsymbol{v}}
        - \int f \frac{ \boldsymbol{F}}{m} \cdot \nabla_{\boldsymbol{v}} g \ d\boldsymbol{v} \Big) \ dE_{\boldsymbol{v}} \\
        & = - \int f \frac{ \boldsymbol{F}}{m} \cdot \nabla_{\boldsymbol{v}} g \ d\boldsymbol{v} \ dE_{\boldsymbol{v}}  \\
        & = - n \langle \frac{ \boldsymbol{F}}{m} \cdot \nabla_{\boldsymbol{v}} g \rangle,
    \end{aligned}
\end{equation}
where $S_{\boldsymbol{v}}$ represent the surface integral in the dimension of $\bfs{v}$.

\item Similarly, we calculate the fourth term to render
\begin{equation}
    \begin{aligned}
        \int g \partial_\omega \Big(f \frac{F_\omega}{I} \Big) \ dE
        & = \int \Big( \int g \partial_\omega \Big(f \frac{F_\omega}{I} \Big) \ d\omega \Big) \ dE_\omega \\
        & = \int \Big( \int g f \frac{F_\omega}{I} \ dS_\omega - \int f \frac{F_\omega}{I} \partial_\omega g \ d\omega \Big) \ dE_\omega \\
        & = - \int f \frac{F_\omega}{I} \partial_\omega g \ dE \\
        & = - n \langle \frac{F_\omega}{I} \partial_\omega g  \rangle,
    \end{aligned}
\end{equation}
where $S_\omega$ represent the surface integral in the dimension of $\omega$.

\item We calculate the last term on the left-hand-side to lead to
\begin{equation}
    \begin{aligned}
        \int g \partial_\theta (\omega f ) \ dE
        = \int \Big( \int \partial_\theta (g \omega f ) \ d\theta \Big) dE_\theta
        = \int \Big( \int (g \omega f ) \ dS_\theta \Big) dE_\theta
        & = 0,
    \end{aligned}
\end{equation}
where $S_\theta$ represent the surface integral in the dimension of $\theta$.

\item Finally, the right-hand-side term is assumed to be conserved in a collision, such as mass, momentum, and energy. Thus,
$
\int g C(f) \ dE = 0
$
as in \eqref{eq:gc0}.

\end{itemize}

Summing up the above equations, we arrive at the generic equation
\begin{equation} \label{eq:gmom}
    \frac{\partial }{\partial t} \big(n \langle g \rangle \big)
    + \nabla_{\boldsymbol{x}} \cdot \big(n \langle g \boldsymbol{v} \rangle \big)
    - n \langle \frac{ \boldsymbol{F}}{m} \cdot \nabla_{\boldsymbol{v}} g \rangle
    - n \langle \frac{F_\omega}{I} \partial_\omega g  \rangle
    = 0.
\end{equation}
To obtain the moment equation of number density, we first set $g=1$, and then
$$
\langle g \rangle = 1, \quad
\langle g \boldsymbol{v} \rangle = \langle \boldsymbol{v} \rangle, \quad
\nabla_{\boldsymbol{v}} g = 0, \quad
\partial_\omega g = 0,
$$
which lead to the continuity equation for the number density
\begin{equation} \label{eq:pde_n}
    \frac{\partial n }{\partial t}
    + \nabla_{\boldsymbol{x}} \cdot \big(n \langle \boldsymbol{v} \rangle \big)
    = 0,
\end{equation}
where $\langle \boldsymbol{v} \rangle$ is the average of floe velocities which are approximated in numerical experiments to be coefficients that are generated as large-scale features from the local DEM models.

Let $g=m$. Using the fact that $n\langle g \rangle = \langle ng \rangle$, the equation \eqref{eq:gmom} reduces to
\begin{equation}
    \frac{\partial \langle \rho \rangle }{\partial t}
    + \nabla_{\boldsymbol{x}} \cdot \big( \langle \boldsymbol{v} \rho \rangle \big)
    = 0,
\end{equation}
where $\rho = \langle ng \rangle$ is a scaled mass density. This is a statement of mass conservation. In the case that floe radius and thickness are fixed, this equation reduces to \eqref{eq:pde_n}.
More importantly, in our setting where the sea ice density $\rho_{ice}$ and sea ice thickness $h$ is fixed,
this equation is equivalent to a conservation statement of concentration $c$:
\begin{equation}
    \frac{\partial \langle c \rangle }{\partial t}
    + \nabla_{\boldsymbol{x}} \cdot \big( \langle \boldsymbol{v} c \rangle \big)
    = 0,
\end{equation}
which is derived from constant scaling.
We will use this formula to update the radii of the local fine-scale DEM models in the next Section.

Now let $g=m \boldsymbol{v}$ and $\boldsymbol{P} = nm \boldsymbol{v} = \rho_{ice} \pi \rho \boldsymbol{v}$. By omitting the scaling constants $\rho_{ice} \pi$, we have
$$
n \langle g \rangle = \langle \rho \boldsymbol{v}  \rangle, \quad
n \langle g \boldsymbol{v} \rangle = \langle \rho \boldsymbol{v} \otimes \boldsymbol{v} \rangle, \quad
\langle \frac{F_\omega}{I} \partial_\omega g  \rangle = 0
$$
and
\begin{equation}
    \langle \frac{ \boldsymbol{F}}{m} \cdot \nabla_{\boldsymbol{v}} g \rangle
    = \langle \boldsymbol{F}  \rangle
    = \langle \boldsymbol{F}^{contact} + \boldsymbol{F}^{drag}  \rangle
    = \langle  \boldsymbol{F}^{drag}  \rangle
    = \langle \widetilde{\alpha}^l\left(\boldsymbol{u}^l_o -\boldsymbol{v}^l\right)\left| \boldsymbol{u}^l_o -\boldsymbol{v}^l\right| \rangle
\end{equation}
where the term associated with contact force vanishes as a result of the conservation of momentum in a collision.
Thus, the equation \eqref{eq:gmom} reduces to
\begin{equation}
    \frac{\partial \langle \rho \boldsymbol{v}  \rangle }{\partial t}
    + \nabla_{\boldsymbol{x}} \cdot \big( \langle \rho \boldsymbol{v} \otimes \boldsymbol{v} \rangle \big)
    = n \langle \widetilde{\alpha}^l\left(\boldsymbol{u}^l_o -\boldsymbol{v}^l\right)\left| \boldsymbol{u}^l_o -\boldsymbol{v}^l\right| \rangle,
\end{equation}
which states that the momentum is conserved with the external drag force under consideration.
\rev{Lastly,} let $g=I\omega$ and $P_\omega = nI\omega$. Then, we have
$$
n \langle g \rangle = \langle P_\omega \rangle, \quad
n \langle g \boldsymbol{v} \rangle = \langle P_\omega \boldsymbol{v} \rangle, \quad
\nabla_{\boldsymbol{v}} g = 0
$$
and
\begin{equation}
    \langle \frac{F_\omega}{I} \partial_\omega g  \rangle
    = \langle F_\omega  \rangle
    = \langle F_\omega^{contact} + F_\omega^{drag} \rangle
    = \langle \widetilde{\beta}^l \left(\nabla\times \boldsymbol{u}_o/2-\omega^l\hat{\boldsymbol{z}}\right)\left|\nabla\times \boldsymbol{u}_o/2-\omega^l\hat{\boldsymbol{z}}\right|  \rangle,
\end{equation}
where the term associated with contact force vanishes as a result of the conservation of angular momentum in a collision.
Thus, the equation \eqref{eq:gmom} reduces to
\begin{equation}
    \frac{\partial \langle P_\omega \rangle }{\partial t}
    + \nabla_{\boldsymbol{x}} \cdot \big( \langle \boldsymbol{v} P_\omega \rangle \big)
    = n \langle \widetilde{\beta}^l \left(\nabla\times \boldsymbol{u}_o/2-\omega^l\hat{\boldsymbol{z}}\right)\left|\nabla\times \boldsymbol{u}_o/2-\omega^l\hat{\boldsymbol{z}}\right|  \rangle,
\end{equation}
which states that the angular momentum is conserved with the drag force taken into consideration.

Higher-order moments may be included for more accuracy. Herein, for simplicity and computational efficiency, we only consider these lower-order moments.
In summary, the moment equations are
\begin{subequations}\label{eq:pde}
\begin{align}
%    \frac{\partial n }{\partial t}
%    + \nabla_{\boldsymbol{x}} \cdot \big(n \langle \boldsymbol{v} \rangle \big)
%    & = 0, \\
    %
    \frac{\partial \langle \rho \rangle }{\partial t}
    + \nabla_{\boldsymbol{x}} \cdot \big( \langle \boldsymbol{v} \rho \rangle \big)
    & = 0, \label{eq:pde_m} \\
    \frac{\partial \langle \rho \boldsymbol{v}  \rangle }{\partial t}
    + \nabla_{\boldsymbol{x}} \cdot \big( \langle \rho \boldsymbol{v} \otimes \boldsymbol{v} \rangle \big)
    & = n \langle \widetilde{\alpha}\left(\boldsymbol{u}_o -\boldsymbol{v}\right)\left| \boldsymbol{u}_o -\boldsymbol{v}\right| \rangle, \label{eq:pde_mv} \\
    \frac{\partial \langle P_\omega \rangle }{\partial t}
    + \nabla_{\boldsymbol{x}} \cdot \big( \langle \boldsymbol{v} P_\omega \rangle \big)
    & = n \langle \widetilde{\beta} \left(\nabla\times \boldsymbol{u}_o/2-\omega\hat{\boldsymbol{z}}\right)\left|\nabla\times \boldsymbol{u}_o/2-\omega\hat{\boldsymbol{z}}\right|  \rangle.   \label{eq:pde_w}
\end{align}
\end{subequations}

We remark that the DEM modeling equations \eqref{eq:dem} are described using the Lagrangian coordinates while the moment equations \eqref{eq:pde} are under the standard Eulerian ones. The unknowns in the above governing equations are averages characterizing large-scale features. \rev{For a 2D simulation domain with a mesh grid, a numerical solution gives the averaged physical quantities in each grid cell where there is a local DEM model for the set of floes inside this cell. 
Hence, the numerical solution gives an approximation of the large scales while the fine-scale features are retained in local DEM models.} This resembles the features of a particle-continuum model. The following section combines the large-scale moment equation model with the fine-scale DEM model\rev{; see Figure \ref{fig:modiag} for a diagram of the particle-continuum model development.}

%As a final remark, 
We also remark that the super-parameterization process resembles the features of the mean-field limit theory for the dynamics of large particle systems, which has been developed in the most detail for molecules or ions and electrons in plasma in \cite{golse2003mean,golse2016dynamics,jabin2014review}. The theory is to derive simplified modeling equations to approximate the system of equations for the dynamics of particles in an averaged sense when the number of particles involved tends to infinity. In fluid mechanics, the particles can also be used to represent vortices in fluid flow simulations \cite{chen1998lattice}, and the corresponding mean field model is the vorticity formulation of the Euler equations of incompressible fluid mechanics. 

%The sea ice particles resemble the features of vortices in fluid dynamics, and the resulting moment equations \eqref{eq:pde} resemble the form of the Navier-Stokes equations.

In relation to sea ice models, the large-scale moment model in \eqref{eq:pde} resembles some aspects of traditional continuum sea ice models \cite{hibler1979dynamic,feltham2008sea}, but it also includes notable differences. One similarity is the use of low-order moments such as mass density and momentum, which is helpful in indicating similarities in basic physical processes and facilitating comparisons of the different modeling approaches. One difference is that \eqref{eq:pde_w} represents the evolution of ice angular momentum, an evolution equation that is not typically included in traditional continuum sea ice models. Another notable difference is that the stress tensor is traditionally specified as a closure assumption and rheology assumption \cite{hibler1979dynamic,feltham2008sea}, whereas here in \eqref{eq:pde_mv} the stress tensor $-\langle \rho \boldsymbol{v} \otimes \boldsymbol{v} \rangle$ is not specified in terms of rheology nor closed in terms of coarse-scale quantities.  Here, the fine-scale floe model and its statistical averages supply the information for the coarse-scale stress tensor $-\langle \rho \boldsymbol{v} \otimes \boldsymbol{v} \rangle$ as part of a multi-scale model.

\section{Particle-continuum coupling at the numerical level} \label{sec:alg}

This section presents the numerical methods for solving the fine-scale and coarse-scale models. In particular, we focus on their coupling and discuss the options and their pros and cons for two-way coupling. We also refer to this coupled model as a multiscale DEM (msDEM) model.

\subsection{Overview}

As the fine-scale model \eqref{eq:dem} is a particle model while the coarse-scale model \eqref{eq:pde} is a continuum one, different numerical methods are developed to solve these two models. Specifically, since the particle model is a set of coupled ordinary differential equations and the physical quantities in time are sampled at the fine time steps, we solve it using the simple forward Euler method with fine time step sizes for accuracy. The coarse-scale model is a time-dependent partial differential equation in 2D. Thus, we need both spatial and temporal discretization schemes. In addition, the \rev{coarse scale} aims to capture the large-scale features of the dynamics. Hence, the coarse-scale model is solved with larger time step sizes and a coarse spatial grid (mesh). With this in mind, we \rev{apply the Lax-Friedrichs scheme \cite{lax1954weak} for spatial discretization.}

\begin{figure}[h!]
    \centering
    \includegraphics[width=0.7\textwidth]{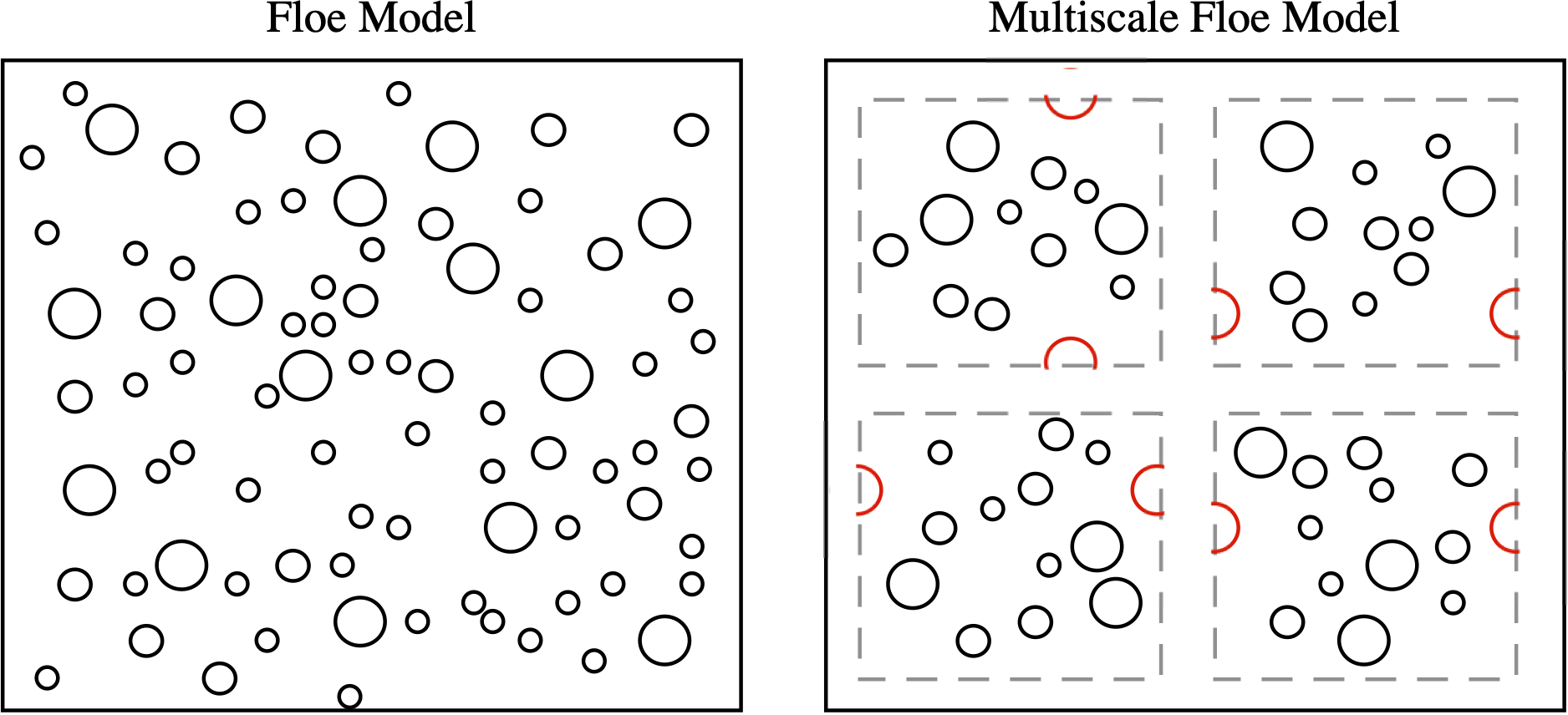}
    \caption{Schematic illustration of a floe model and a multiscale floe model. In the multiscale floe model, only four coarse cells (in a $2\times 2$ arrangement) are illustrated in this schematic diagram. Within each coarse cell is a floe model.  The red floes illustrate that doubly-periodic boundary conditions are used for the floe model within each coarse cell.}
    \label{fig:multiscale-schematic}
\end{figure}

To couple the models, it is essential to distinguish the different frames of reference used in both models. The particle model uses a Lagrangian frame of reference, while the continuum model uses the Eulerian coordinates. Figure \ref{fig:multiscale-schematic} shows a schematic illustration of the setup of the 2D spatial grid for the multiscale model. Within each coarse grid cell, we solve one fine-scale DEM particle model. The main idea of the two-way coupling is as follows. In each \rev{large-scale grid cell}, we run a DEM particle model for a set of floes with double-periodic boundary conditions. \rev{Herein, we apply the double-periodic boundary conditions for fully parallel computing (as in typical superparameterization for the atmosphere; e.g., \cite{grabowski1999crcp, randall2003breaking, majda2009simple, stechmann2014multiscale,  grooms2013efficient, majda2014new, lee2017stochastic, dolaptchiev2013multiscale, jansson2019regional}). With these conditions, one does not need to transfer floes with the DEM models in neighbouring processors in a parallel computing environment. The communications among processors are established by solving the continuum model for large-scale physical quantities and updating the local DEM models through model coupling. This significantly reduces the communication time among processors. In the numerical section, we demonstrate the convergence behaviour of using this periodic boundary condition for the multiscale model to the full DEM model.  } 
The floes' sizes, velocities, and angular velocities are adjusted according to the cell averages that are solutions of the coarse-scale model. 
\rev{We initialize the system so that each local DEM has a statistically significant number of particles.  }
We solve the moment equations for the large-scale physical quantities on the overall domain with a coarse grid while the right-hand side forces in \eqref{eq:pde} are inputs processed from the fine-scale particle solutions. The large-scale solutions provide gain or loss on total mass (characterized by radius as thickness and density are constants), total momentum, and total angular momentum for each coarse-grid cell where the DEM lives in. The fine-scale DEM model provides the coefficients and external forcing (ocean drag) for the moment equations, and these coefficients are given by mainly the averaged quantities from the DEM models, such as averaged velocity for the mass density equation in \eqref{eq:pde}.

\begin{figure}[h!]
\centering
	\includegraphics[width=6in]{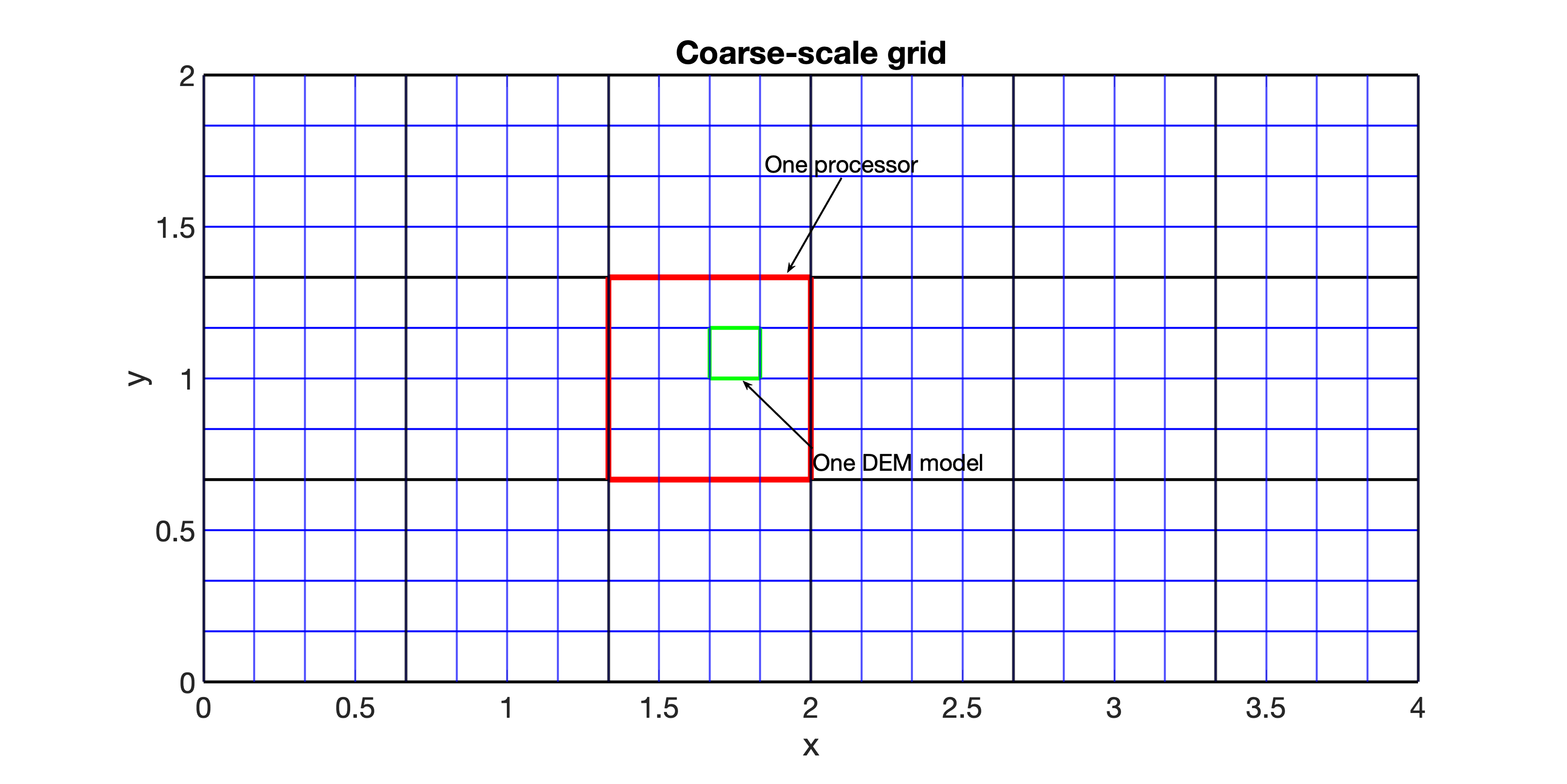}
	\caption{An example of HPC processor and coarse-scale grid setting. There are $6\times3=18$ processors (black and bold lines) and $24\times12$ coarse-grid cells (blue and thin lines). Each coarse-grid cell contains one DEM particle model (a representative one is marked in green).}
	\label{fig:cgrid}
\end{figure}

\begin{algorithm}[h!]
\caption{Coarse- and fine-scale model coupling}
\label{alg:cp}
\begin{algorithmic}
\State{Initialize a system with $N\geq1$ floes.}
\State{Initialize time-marching parameters $dt, N_t, N_0, N_1, n_t$.}
\State{Run the local DEM models with the fine time step size $dt$ to the \rev{coarse time step $dT$}. }
\State{Collect the time series data over time steps $j=1,2,\cdots, N_0-1$. }
\State{Initialized the coarse-scale model with the DEM data.}
\For{a fine time step marching $j=N_0, \cdots, N_t$ with final time $T = N_t dt$}

\State{Accumulate local DEM data}

\If{$\mod(j, N_0)=0$}

\State{Process the accumulated DEM data to obtain input for the coarse-scale model.}
\State{Solve the large-scale model to time $(j+N_0)dt$ with $N_1$ substeps such that $dT=N_1 \cdot (N_0 dt/N_1)$.}
\Else

\If{$j>N_0$ and $\mod(\mod(j, N_0), n_t)=0)$}
    \State{update the local DEM floes gradually using the large-scale model solutions}
\EndIf

\EndIf
\EndFor

%\While{Number of coarse time step}
%\EndWhile
\State{Save data for postprocessing.}
\end{algorithmic}
\end{algorithm}

Figure \ref{fig:cgrid} shows an example of a coarse-scale grid. We consider non-dimensional units in this paper. The computational global domain of Figure \ref{fig:cgrid} is $D = [0,4]\times[0,2]$. The coarse grid has $24\times 12 = 288$ cells, each containing one fine-scale DEM particle model. For high-performance computing (HPC) purposes, we solve both the coarse and fine-scale models in a parallel environment with many processors.  In the example shown in Figure \ref{fig:cgrid}, there are $18$ processors that are uniformly distributed on the simulation domain. The HPC grid is of size $6\times3$. In each processor, there are $4\times 4$ grid cells meaning $4\times 4 = 16$ DEM particle models that are solved independently with double-periodic boundary conditions. More details are given in the next subsection below, with implementation details on how the coupling is realized on HPC clusters. \rev{Lastly, the two-way coupling is} summarized in Algorithm \ref{alg:cp}.

\subsection{Coarse-scale moment equations}

The equations \eqref{eq:pde} are first-order hyperbolic PDEs and are not supplied with boundary conditions. Due to the divergence (first-order partial derivative) structure, spatial discretization schemes such as finite element methods often suffer stability issues (challenging to design schemes that satisfy the Ladyzhenskaya–Babuška–Brezzi (LBB) inf-sup condition). This generally leads to difficulties in applying implicit time-marching schemes for solving the modeling equations \eqref{eq:pde}. \rev{Moreover, implicit time-marching schemes require inverting matrix systems in each time step, which can be time-consuming. Implicit time-marching schemes have robust stability and enable larger time steps. However, herein we use small time step sizes to ensure numerical accuracy in the temporal dimension. } Thus, we herein \rev{use Lax-Friedrichs scheme \cite{lax1954weak} that is widely used for solving hyperbolic PDEs.}

%explicit time-marching schemes, particularly the two-step Adams-Bashforth method \cite{hairer1993solving}. For spatial discretization, we use the Lax-Friedrichs scheme \cite{lax1954weak} that is widely used for solving hyperbolic PDEs.

We denote the fine-scale and coarse-scale time-marching steps by $dt$ and $dT$. To resolve the floe collisions in DEM fine-scale model and obtain numerical accuracy, small time steps are required, and we set in the default $dt=10^{-4}.$ We set it based on numerical stability and accuracy for the coarse-time step size. For a coarse mesh element size $dX$ and approximate velocity $v$, the CFL condition requires roughly
    \begin{equation}
        \Big| \frac{v dT}{dX} \Big| < C_{max},
    \end{equation}
where $C_{max}$ is the Courant number. The default value is $C_{max}=1$ for an explicit time-marching scheme. Since $dT >0, dX>0,$ this leads to
    \begin{equation}
        dT < \frac{C_{max} dX}{ |v| }.
    \end{equation}
\rev{For example, on a global domain $[0,4] \times [0,2]$ with a grid of size $64\times32$ and maximum velocity $v = 1.0$, the time-step size for numerical stability is $dT < \sqrt{2}/16 \approx 0.09$.}
For better numerical accuracy, we set in default $dT=0.01.$ We denote $N_0 = \frac{dT}{dt}.$

To initialize the solutions of moment equations \eqref{eq:pde}, we let the local DEM models run for time steps $j=1,\cdots, N_0$ to get DEM statistical data. This includes the scaled mass density, averaged floe velocities, and averaged floe angular velocities. To couple the DEM fine-scale particle model, we perform the following steps.

\begin{enumerate}
    \item We assume that there are $n$ floes in each coarse cell. It is given as the total number of particles in this coarse cell and it is fixed in the overall simulation. There is one value for each physical quantity, $\rho, \rho \boldsymbol{v}, P_\omega$, in one cell. These values are governed by the coarse-scale model \eqref{eq:pde}. They represent averaged values of all the $n$ floes in that cell and we use them as a base for updating the local DEMs.

    %The mass density or floe concentration is updated by updating floe radii.

    \item To solve \eqref{eq:pde_m}, we treat $\langle \boldsymbol{v} \rangle$ as a known quantity, which is an average of the floe velocities in the coarse cell. This is an average over fine-scale time steps obtained from the local DEM models.

    \item To solve \eqref{eq:pde_mv} and \eqref{eq:pde_w}, the unknowns are $\langle \rho \boldsymbol{v} \rangle$ and $\langle P_\omega \rangle$, respectively. The other terms involving velocities and the right-hand side terms are treated as known values approximated from the fine-scale DEMs.
    In particular, for the momentum equation \eqref{eq:pde_mv}, we solve it using the idea adopted in the derivation of Oseen equations (for approximating the Navier-Stokes equation). We treat the part of the quadratic advection term as a known quantity generated from the local DEMs. That is, we approximate the term
    \begin{equation}
    \nabla_{\boldsymbol{x}} \cdot \big( \langle \rho \boldsymbol{v} \otimes \boldsymbol{v} \rangle \big)
    \approx
    \langle \boldsymbol{v} \rangle \cdot \nabla_{\boldsymbol{x}} \big( \langle \rho \boldsymbol{v} \rangle \big).
    \end{equation}

\end{enumerate}

\subsection{Small-scale DEM}

To solve the fine-scale DEM particle model \eqref{eq:dem} in a coupled way, we use the large-scale solutions to gradually update the local DEMs' physical quantities gradually over an equally spaced sequence of time nodes within one coarse time step. The DEM model is solved for floe positions, velocities, and angular velocities at all fine time step $jdt, j=1,2,\cdots, N_t$. The coarse-scale model \eqref{eq:pde} is solved for cell averages of floe (scaled) mass density (equivalent to concentration $c$), velocities, and angular velocities at all coarse time step $kdT, k =1,2, \cdots, \frac{N_t}{N_0}$ where $dT = N_0 dt$. Each coarse time step has $N_0$ fine time steps. These $N_0$ fine time steps are further divided into $n_t$ steps for gradual updates on local DEMs. In default, we set \rev{$N_0=100$} and $n_t=10$, which means we gradually update the local DEM floes over $n_t=10$ fine time steps in each coarse time step.

To sum up, with the coarse-scale solutions, we perform the following steps to solve the local DEM particle models while gradually updating the averages of floe physical quantities.

\begin{itemize}
    \item We integrate the scaled mass density $\rho$ over the local DEM domain, denoted as $D_{DEM}$, to get the averaged floe concentration in the domain. That is, the number of particles in the DEM domain is
    \begin{equation}
        \bar{r} = \frac{\int_{D_{DEM}} \rho \ d\boldsymbol{x}}{\rho_{ice} \pi \int_{D_{DEM}} 1 \ d\boldsymbol{x}}.
    \end{equation}
    The target average floe radius of all the floe particles in this coarse grid cell in the next time-marching time step is thus $\bar{r}$. We consider a coarse time step from $t = kdT$ to $t = (k+1)dT$ with $n_t$ sub-steps. To achieve the target $\bar{r}$ as an average (coarse-scale value), we proportionally distribute the increment $\delta r = \bar{r} - \bar{r}^0$ to each floe according to their radii and update the system to reach this target gradually over $n_t$ steps. $\bar{r}^0$ is the current averaged radius of all the floes in that cell. Thus, we update each floe as follows
    \begin{equation}
        r^{l,new} = r^l + n(\bar{r} - \bar{r}^0 ) \frac{r^l}{\sum_{j=1}^{n} r^j} \frac{j n_t}{N_0}, \quad l = 1,2,\cdots,n, j=1,2,\cdots,\frac{N_0}{n_t},
    \end{equation}
    where $n$ is the total number of floes in this cell. With this update, we arrive at the expected averaged (large-scale) radius at each coarse time step while keeping the variations of individual floes in the fine scale.

    \item Similarly, to update the averaged velocities in each coarse-grid cell, we consider a coarse time step from $t = kdT$ to $t = (k+1)dT$. From the large-scale solution $m \boldsymbol{v} = \boldsymbol{P} = (P_x, P_y)$, it gives the target averaged momentum overall all floes in the coarse grid cell. To achieve this target as an average (coarse-scale value), we proportionally distribute the increment $\delta \bfs{P} = \boldsymbol{P} - \boldsymbol{P}^0$ to each floe according to their masses and update the system to reach this target gradually over $n_t$ steps. $\boldsymbol{P}^0$ is the current averaged momentum of all the floes in that cell. Thus, we update each floe as follows
    \begin{equation}
        \boldsymbol{v}^{l,new} = \boldsymbol{v}^l + \frac{n(\boldsymbol{P} - \boldsymbol{P}^0) }{m^l} \frac{m^l}{\sum_{j=1}^{n} m^j} \frac{j n_t}{N_0}, \quad l = 1,2,\cdots,n, j=1,2,\cdots,\frac{N_0}{n_t},
    \end{equation}
    where $n$ is the total number of floes in this cell. The angular velocities are updated in a similar fashion; that is
    \begin{equation}
        \omega^{l,new} = \omega^l + \frac{n(P_\omega - P_\omega^0) }{I^l} \frac{I^l}{\sum_{j=1}^{n} I^j} \frac{j n_t}{N_0}, \quad l = 1,2,\cdots,n, j=1,2,\cdots,\frac{N_0}{n_t},
    \end{equation}
    where $P_\omega$ and $P_\omega^0$ are the targets and current average angular momentum.
    \rev{With these updates for each of the particles on the local DEMs, we arrive} at the expected averaged (large-scale) velocities at each coarse time step while keeping the variations of individual floes in the fine scales. This allows for conserving momentum.

\end{itemize}

\section{Numerical experiments} \label{sec:num}
In this section, we perform several numerical tests.
%In particular, we start with a simple model scenario where all floes are traveling in one direction without floe contact.
For the error study, one simulation from the full DEM model on the entire domain is regarded as the truth.

\subsection{A simple simulation scenario} \label{sec:ex1}
We first consider a simple scenario with a non-dimensional domain $\Omega=[0,4] \times [0,2]$ and a constant ocean velocity field $\bfs{v} = (0.3, 0)$. We discretize the domain into a mesh with $480\times 240$ uniform grid cells (squares, elements). In each cell, we place in the center a floe with a radius given as
\begin{equation}
    r = r_c (0.2 + 0.8\sin(0.25\pi x)),
\end{equation}
where $r_c = \frac{1}{240}$ such that $2r_c$ is the cell length. Hence, there are $480\times 240 = 115200$ floes in total. The initial velocities of floes are set to be \rev{$\bfs{u}_o = (0.3, 0)$}. We apply periodic boundary conditions. Thus, in this simple scenario, the floes travel to the right periodically without floe-to-floe contact. The angular velocities are initialized to zeros and remain zeros throughout the simulation, \rev{both in DEM (benchmark full particle model) and msDEM models (the proposed particle-continuum multiscale model).} For this simulation scenario, we focus on the floe concentration.

 \begin{figure}[h!]
\centering
	\includegraphics[width=15cm]{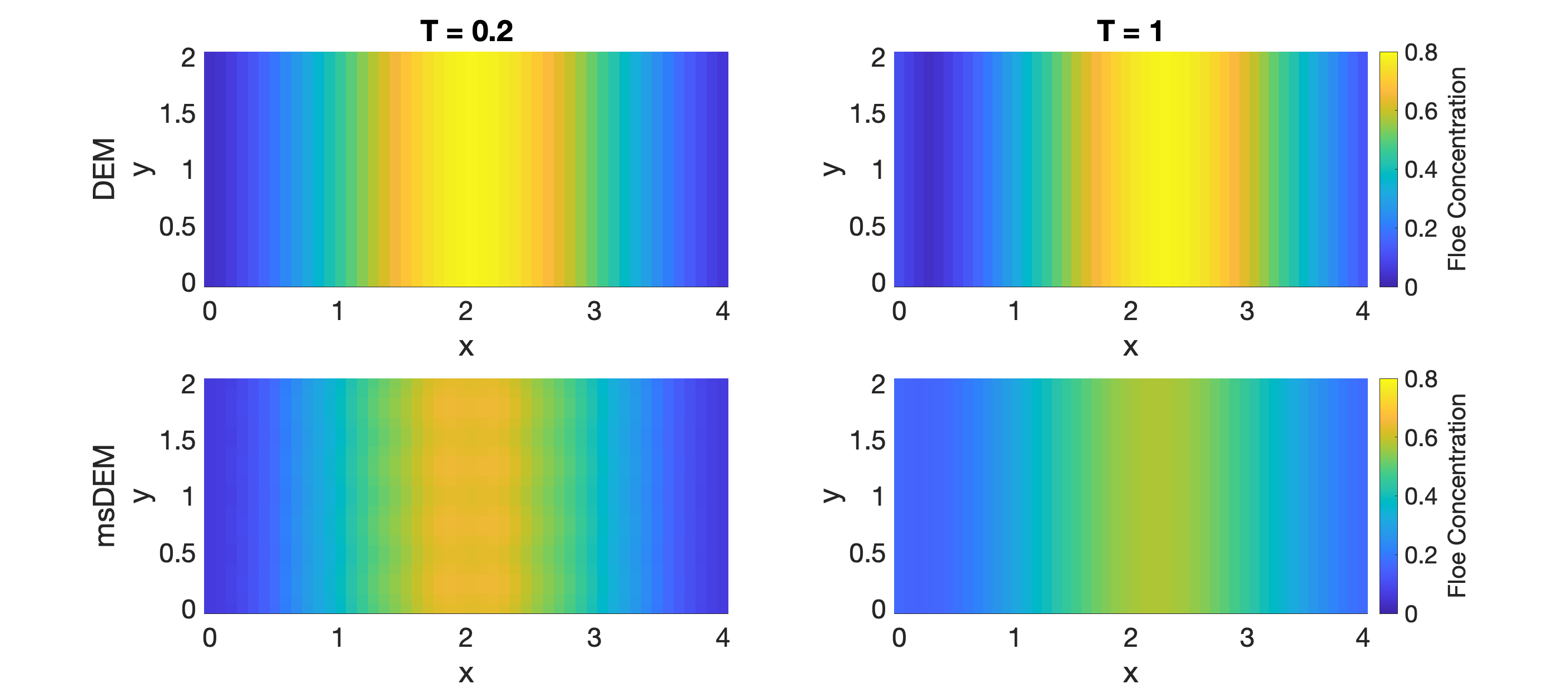}
	\caption{Scenario 4.1. Floe concentration of DEM and msDEM models at $T=0.2$ and $T=1$. The grid is \rev{uniform and of size $48\times 24$ on the domain $\Omega=[0,4] \times [0,2]$}.}
	\label{fig:floec}
\end{figure}

Figure \ref{fig:floec} shows the floe concentrations of DEM and msDEM models over a coarse grid of size $48\times 24$ at time $T=0.2$ and $T=1$.
With the setting described above in mind,
each grid cell is of size $\frac{1}{12} \times \frac{1}{12}$ and contains 100 floes.
The concentration is calculated as a fraction of the total area of these 100 floes over the area of the grid cell that is $\frac{1}{144}.$
\rev{The numbers in the colour bars indicate the concentration; same for other figures in the simulations below.}
The msDEM captures the overall dynamic features of the DEM solution, especially for a relatively shorter simulation time. More importantly, the msDEM model is fully parallelized, allowing a significant reduction in the computational cost. In our simulation, we observe that the msDEM model with 32 processors is about six times faster than the DEM model with four processors, leading to a parallel efficiency of approximately $75\%$. We observe similar parallel efficiency when using other settings of processors in the MPI (message passing interface) parallelization and other simulation scenarios discussed below. The parallel efficiency depends on a few aspects: (a) the initial setting of the floes, (b) the ocean-driven forces (which affect the computationally-demanding floe collision calculation), and (c) the parallel setting. This requires an extensive study and is not the focus of this paper. Therefore, we briefly state our simulation observation without rigorous numerical validation.

\begin{figure}[h!]
\centering
	\includegraphics[width=12cm]{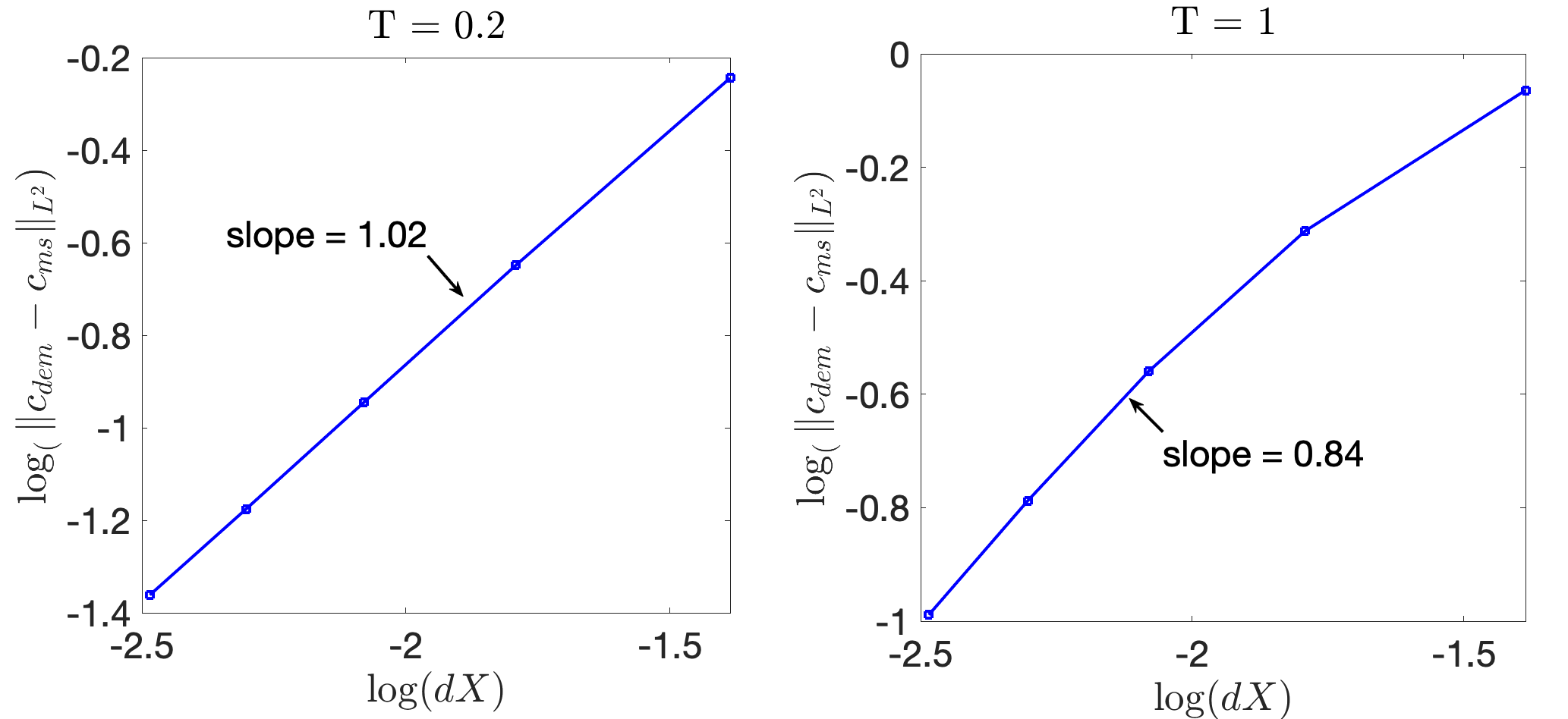}
	\caption{Scenario 4.1. $L^2$ errors with respect to coarse grid size $dX$ for concentration $c$. }
	\label{fig:ex1cerrT}
\end{figure}

Figure \ref{fig:ex1cerrT} shows the $L^2$ errors (in logarithmic scale) of the floe concentration with respect to the coarse-scale model grid size at a few different times.
For a real-valued function $u$ defined on $\Omega$, the $L^2$-norm is defined as $\| u \|_{L^2} = \sqrt{\int_\Omega u^2 \ dx}$.
In the discrete setting, the integral is calculated as a summation over each grid cell.
Therein, $c_{dem}$ refers to the concentration field (with respect to the coarse-scale grid) of the DEM model, while $c_{ms}$ refers to that of the msDEM model.
We observe that for a short time $T$, the convergence rate is of order one, which is optimal as we applied the Lax-Friedrichs scheme for the spatial discretization \cite{lax1954weak}.
\rev{
The numerical solution for the continuum model on coarse grids leads to large errors. 
When using finer grids, we need a larger number of particles so that each cell has a statistically significant number of floes (this affects the accuracy of the coupling as averages of physical quantities in cells are utilized).
}
For longer times, we observe larger errors and diffusion behavior on the msDEM model. This is due to the following three reasons. First, the DEM model itself has time-step accumulating errors, which leads to larger errors when evaluating the concentration in each coarse grid cell. Second, the coarse-scale model takes the statistic moments of the fine-scale DEM models, which averages the floe concentrations. Third, the numerical scheme \rev{(Lax-Friedrichs scheme)} used for solving the coarse-scale model \eqref{eq:pde} has a dissipation behavior. In reality, the DEM particle models have usually incorporated data assimilation (DA) techniques for higher accuracy in longer-time simulations. One would expect higher accuracy in the DA-assisted msDEM models, but its development is out of the scope of this paper, and it is subject to future work.
Lastly, we observe that floes do not contact each other in the msDEM simulation. This is consistent with the setup that no floe collisions exist in this particular DEM simulation scenario.

\subsection{Gathering and scattering floes}

In the previous subsection, we demonstrated the convergent behavior of the floe concentration of the super-parameterized DEM model. We now vary the ocean velocity and study the performance of msDEM on the large-scale floe velocities. Specifically, we adopt the simulation settings in subsection \ref{sec:ex1} and only revise the ocean velocity to be \rev{$\bfs{u}_o = (0.3-0.1\cos(\pi x), 0)$.} The ocean velocity is one-dimensional and induces floes to gather and scatter over time. Moreover, the ocean, in this case, is \rev{compressible. Angular} velocities remain zeros because (a) the curl of the ocean velocity leads to a zero vector, and (b) the tangential force remains zero due to orthogonality and zero initial angular velocities. We observe zero angular velocities in our simulations in both DEM and msDEM models.

\begin{figure}[h!]
\centering
	\includegraphics[width=15cm]{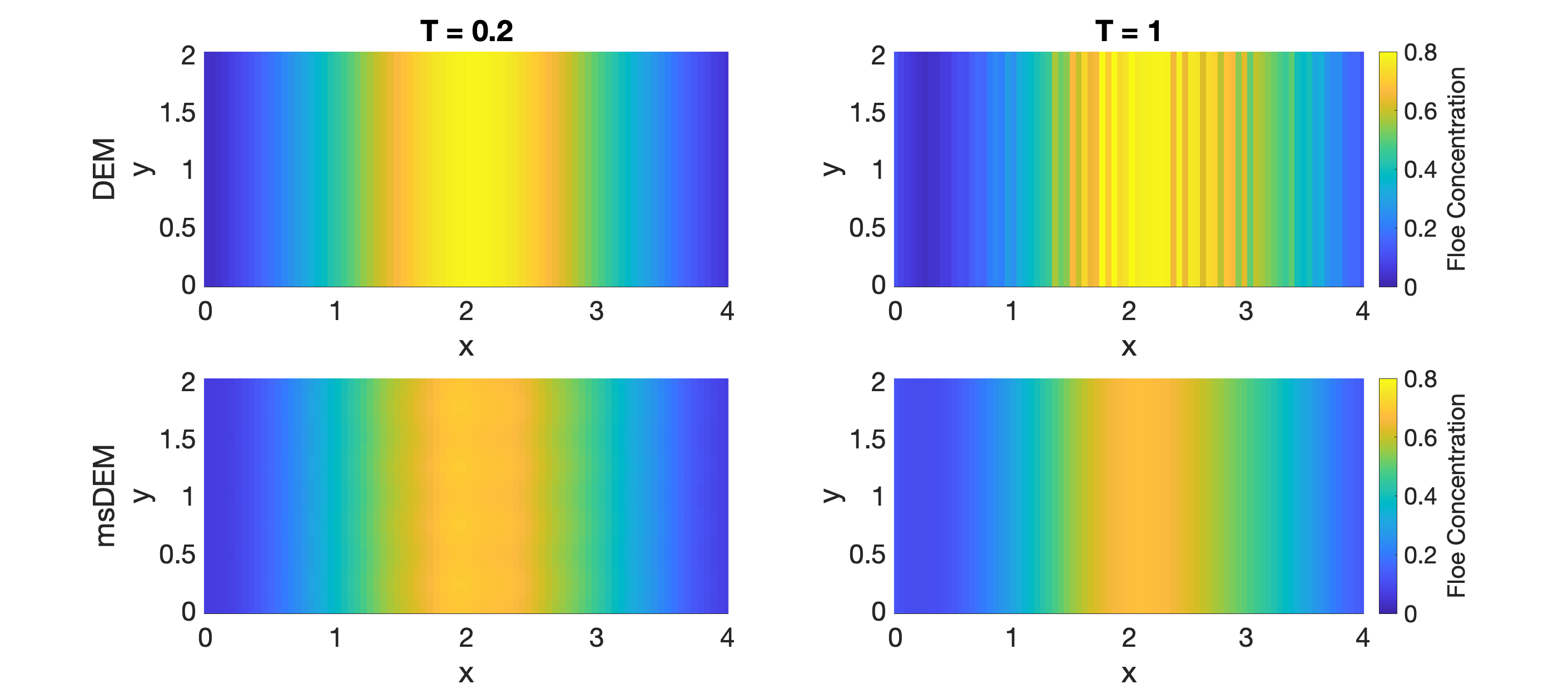}
	\caption{Scenario 4.2. Floe concentration of DEM and msDEM models at $T=0.2$ and $T=1$. The grid is \rev{uniform and of size $48\times 24$ on the domain $\Omega=[0,4] \times [0,2]$}.}
	\label{fig:floecex2}
\end{figure}

\begin{figure}[h!]
\centering
	\includegraphics[width=12cm]{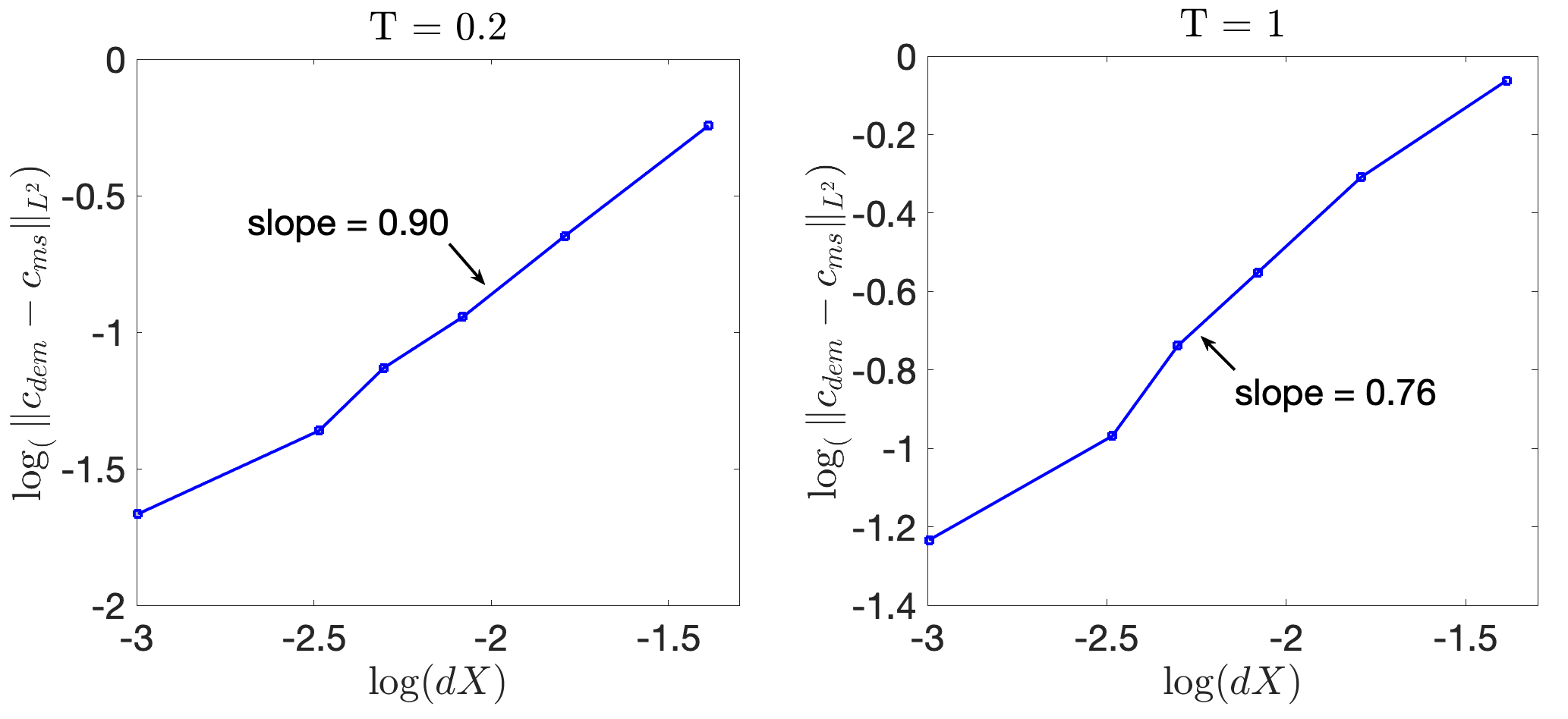}
	\caption{Scenario 4.2. $L^2$ errors with respect to coarse grid size $dX$ for concentration $c$. }
	\label{fig:ex2cerr}
\end{figure}

Figure \ref{fig:floecex2} shows the floe concentration while Figure \ref{fig:ex2cerr} shows the $L^2$ error convergence behavior.
The results are collected following the same procedure as in simulation Scenario 4.1.
We observe that the error converges quickly as the coarse grid gets finer. The plot for the concentration of the DEM simulation at $T=1$ is not smooth because the floes at a particular time can be at locations of one of the grid cells. When floes are near the boundary, depending on the dynamics at a time instance, they possibly belong to one of the neighboring grid cells.
For a short time, the convergence orders are close to one and smaller than those observed in simulation Scenario 4.1. The contributing aspects of smaller convergence orders are: (a) there are collisions in this scenario while there are no collisions in Scenario 4.1, which leads to a more complex dynamical system for msDEM simulation to capture; (b) the floe collision in DEM leads to more errors when evaluating the concentration (see the non-smoothness of the concentration in the top-right plot of Figure \ref{fig:floecex2}).

\subsection{Floe dynamics with a mild compressible ocean current}

\begin{figure}[h!]
\centering
	\includegraphics[width=8cm]{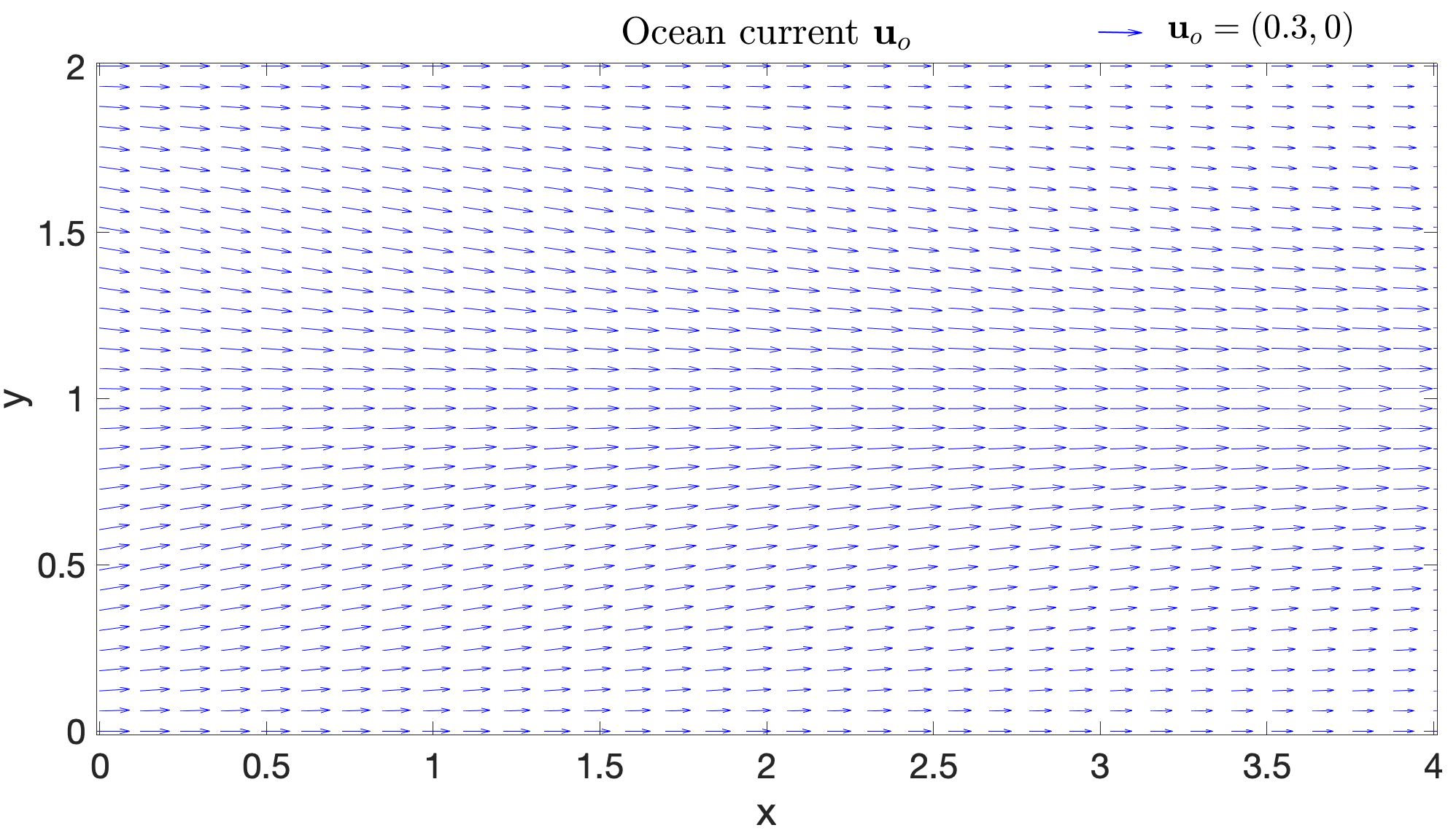}
	\includegraphics[width=8cm]{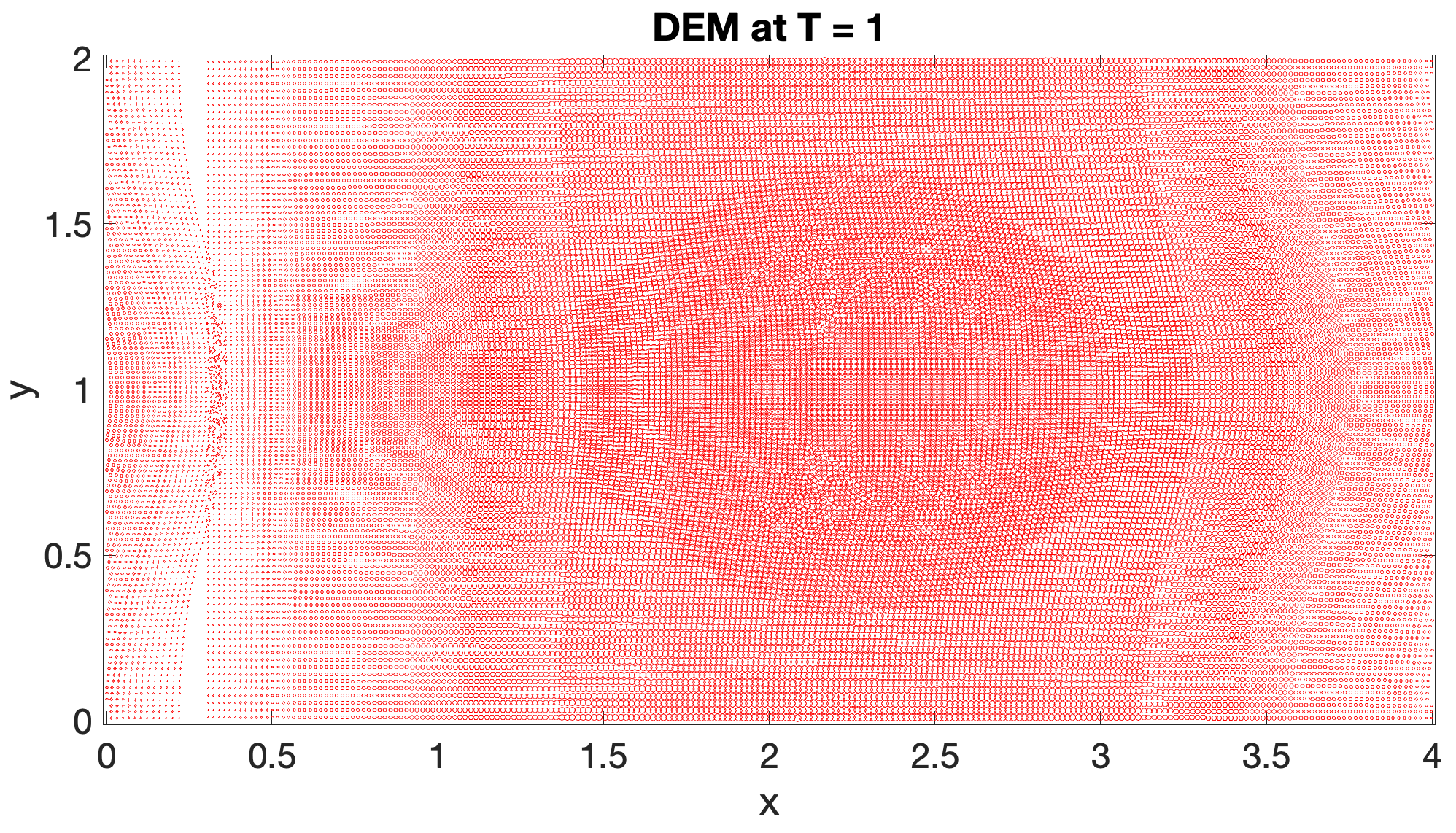}
	\caption{Scenario 4.3. A mild compressible ocean velocity field and the DEM simulation at $T=1$ with 28800 floes.}
	\label{fig:ex3}
\end{figure}

\begin{figure}[h!]
\centering
	\includegraphics[width=15cm]{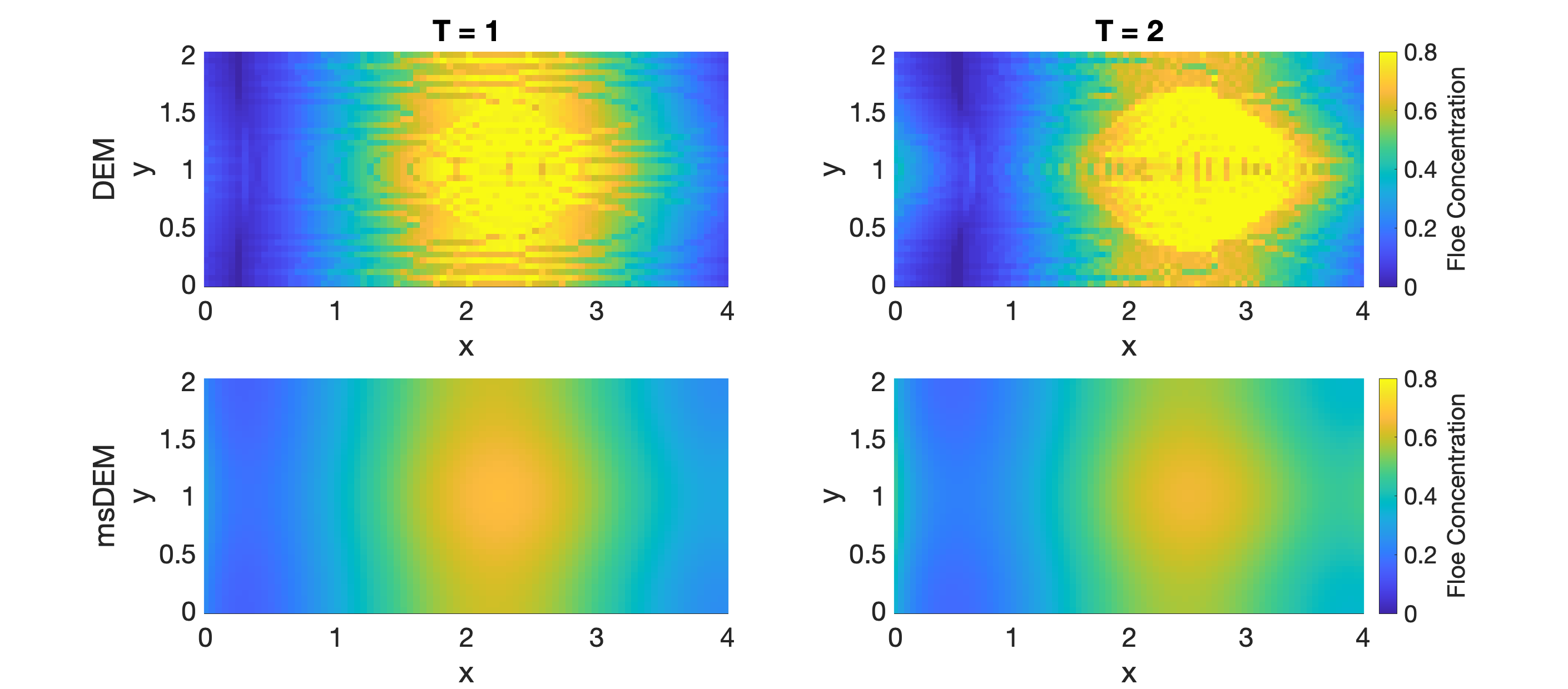}
	\caption{Scenario 4.3. Floe concentration of DEM and msDEM models at $T=1$ and $T=2$. The grid is \rev{uniform and of size $48\times 24$ on the domain $\Omega=[0,4] \times [0,2]$}.}
	\label{fig:ex3c}
\end{figure}

Now, we consider a scenario with a mild compressible ocean current. The ocean velocity field is given as
\rev{
$$
\bfs{u}_o = (0.3-0.1\sin(0.1\pi x) \cos(\pi y), 0.05\cos(0.1\pi x) \sin(\pi y)).
$$}
The left panel of Figure \ref{fig:ex3} shows the velocity field. 
%We observe that the current moves to the center of the horizontal direction, mimicking more realistic ocean currents \rev{(compared to the currents above)}. 
The divergence $\nabla \cdot \bfs{v} = 0.04\pi \cos(0.1\pi x) \cos(\pi y)$ is of relatively small scale, thus a mild compressible ocean current. The initial setting is the same as in scenario 4.1. The right plot of Figure \ref{fig:ex3} shows floes of the DEM simulation at time $T=1$. Herein, we use $240\times120 = 28800$ floes to plot floes (in the simulation, we used 115200 floes as \rev{before. The figure is not readable if using 11520 floes; we show a case with a smaller number of floes to show the state of the floe distribution}). The top-left plot of Figure \ref{fig:ex3c} shows the floe concentration over the coarse mesh grid $48\times 24$.

\begin{figure}[h!]
\centering
	\includegraphics[width=12cm]{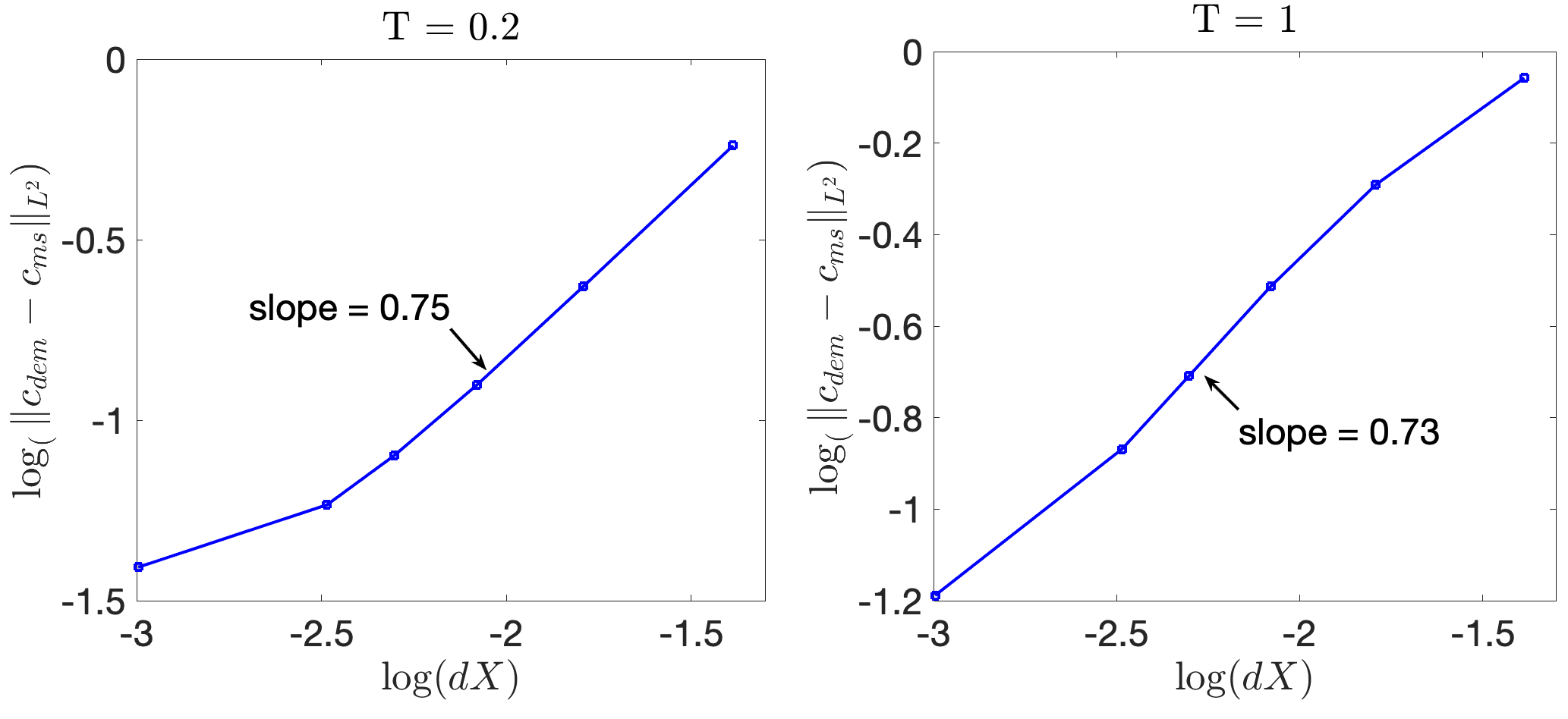}
	\caption{Scenario 4.3. $L^2$ errors with respect to coarse grid size $dX$ for floe velocities $v_x$. }
	\label{fig:ex3err}
\end{figure}

Figure \ref{fig:ex3c} shows the comparison of floe concentration of the DEM and msDEM simulations. We observe that the msDEM overall captures the DEM simulated features well, although the floe concentration is not perfectly smooth due to the floe number counting in each grid cell for floes near the cell boundaries. Figure \ref{fig:ex3err} shows the behavior of the $L^2$ errors of the floe concentration. \rev{Despite that the error is affected more by the non-smoothness of DEM concentrations, a similar error convergent behaviour with a rate around 0.75. }
The convergence orders are getting smaller than the simpler simulation Scenarios 4.1 and 4.2. This is expected as the ocean current introduces non-zero angular velocity components, leading to a more complex dynamic system for msDEM models to simulate. One possible way to increase the msDEM accuracy (hence the convergence rates) is to run simulations with significantly more floes, in which case the concentration field of the DEM models would be much smoother. We keep using 115200 floes in our simulation as a significant increase in the number of floes would require a substantial increase of computational resources, both in storage and time.

\subsection{Floes smashing to a wall}

In this last simulation scenario, we consider the case where floes are smashing into a fixed wall. The initial setting remains the same as in scenario 4.1, and the ocean velocity field is fixed as \rev{$\bfs{u}_o = (0.3, 0)$}. For the $y$-dimension, we apply periodic boundary conditions. For the $x$-dimension, we use a homogeneous (zero) Dirichlet boundary condition on the floe concentration. When a floe moves towards the right boundary, it \rev{gets accumulated}. Thus, the concentration increases near the boundary.
\begin{figure}[h!]
\hspace{-1.8cm}
	\includegraphics[width=20cm]{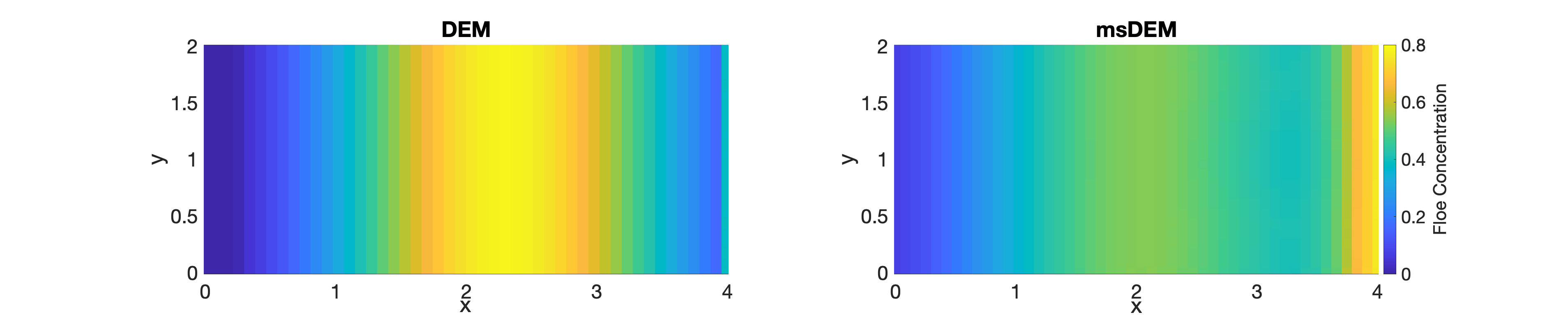}
	\caption{Scenario 4.4. Floe concentration of DEM and msDEM models at $T=1$. The grid is \rev{uniform and of size $48\times 24$ on the domain $\Omega=[0,4] \times [0,2]$}.}
	\label{fig:ex4c}
\end{figure}
\begin{figure}[h!]
\centering
	\includegraphics[width=15cm]{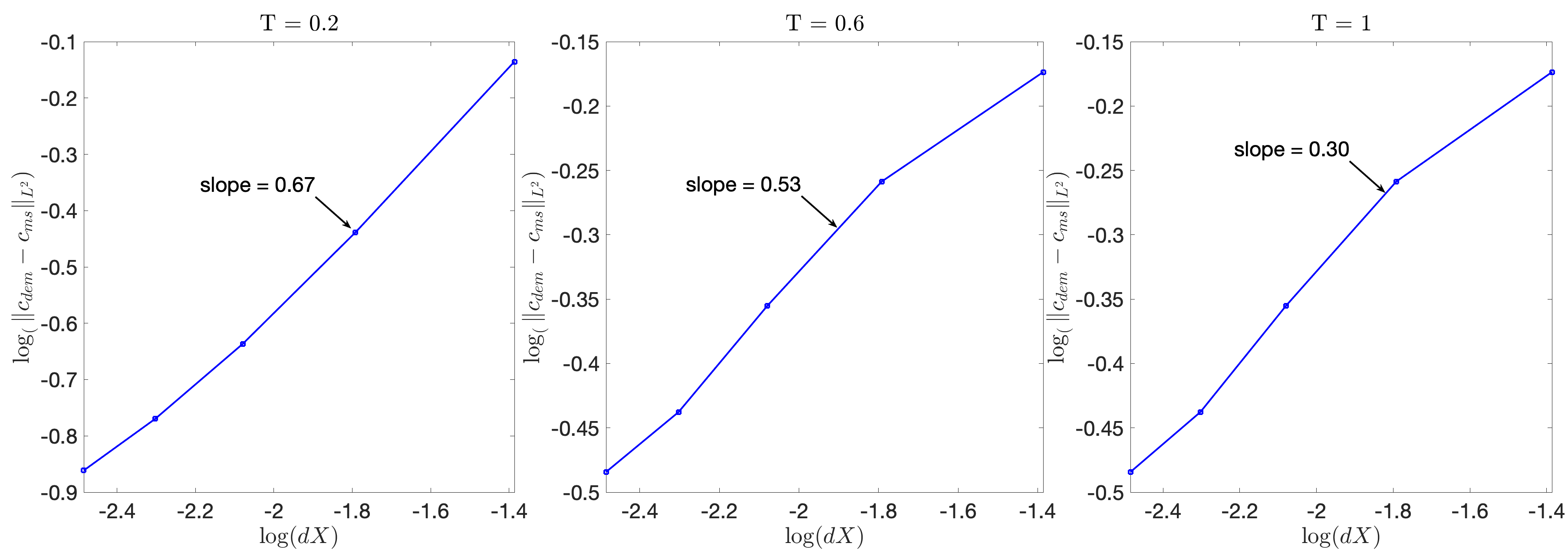}
	\caption{Scenario 4.4. $L^2$ errors with respect to coarse grid size $dX$ for floe velocities $v_x$. }
	\label{fig:ex4err}
\end{figure}
Figure \ref{fig:ex4c} shows the comparison of the floe concentration of the DEM and msDEM simulations at $T=1$.
The msDEM simulation captures the main feature of floes smashing into a wall. Although it accumulates a bit faster than the DEM near the right boundary due to mesh resolution, we expect more accuracy when using finer coarse grids for the large-scale continuum part of the msDEM model.
Figure \ref{fig:ex4err} shows the errors and their convergence rates at $T=0.2,0.6,1$. Due to the boundary layer, the errors are more significant, and convergence rates are lower than in the scenarios considered above.
This scenario is a boundary layer problem, which has intrinsic difficulty for accurate particle-based simulations. The floes get jammed near the right boundary leading to more floe-floe collisions. To better capture the dynamic behavior using the particle-continuum coupled msDEM model, one possible idea is to adopt an adaptive (coarse) mesh/grid which better captures the physical quantities near the boundary layer. However, the adaptive mesh requires a different parallel computing environment, subject to future work. Another point is that once the floes accumulate near the boundary, they become stationary, which could be considered a new boundary with modified physical properties such as Young's moduli. In such a situation, the computation cost could be saved significantly as no more computations are needed for these stationary floes when accumulating near the boundary.

With all these simulation scenarios in mind, we briefly summarize the performance of the proposed particle-continuum msDEM model. All four simulation scenarios demonstrate a feasible and, to some extent, accurate multiscale model to capture the main dynamical feature of a system consisting of a relatively large number of individual particles/floes. The DEM model would be the preferable option for a system with a small number of particles ($\lesssim 10^4$). In contrast, the Boltzmann method based on the statistical distribution of particles is preferred for a system with a considerably large number of particles ($\gtrsim 10^9$). This msDEM model provides an alternative and efficient method to simulate the dynamical behavior of a system with any number of particles in between. 
In a MIZ, depending on how big the region we consider and the scale we want to capture,
the typical number of sea ice floes can be of scale $10^3$ to $10^{12}$ \rev{(estimated based on floe size distributions \cite{horvat2015prognostic,horvat2017evolution,denton2022characterizing})}. The proposed msDEM model inherits the advantages both from the particle model and the continuum model. On the other hand, this multiscale paradigm provides a naturally scale-separately model for data assimilation. In particular, the coarse-scale part of the msDEM model can be used to assimilate the large ocean or atmosphere modes, while the particle-based fine-scale part can be used to assimilate the small ocean or atmosphere modes.

\section{Concluding remarks} \label{sec:con}

The accurate multiscale characterization of the particle simulations remains a significant challenge for sea ice DEM particle simulations. The main difficulty lies in the particle description of the fine-scale dynamics. At the same time, the large scale looks into the averaged dynamical behavior of groups of particles, best described by a continuum model. This work initializes this line of work and adopts the super-parameterization to generate a continuum model that governs the large-scale features of floe particles. In particular, we focus on floe concentration. 
The fine- and coarse-scale models are coupled through averaged physics quantities over each coarse-grid cell. We studied the convergence behavior of the proposed super-parameterization numerically. As far as we know, this is among the first studies that (1) uses a super-parameterization technique to couple the coarse-scale concentration with fine-scale sizes of individual floe particles; and (2) studies the convergence behavior of a super-parameterization technique.
\rev{This multiscale framework is developed for a simplified DEM particle model with assumptions such as no ridging and no atmospheric winds. 
In particular, ridging affects the floe thickness distribution and 
it also leads to out-of-plane motions. 
The generalisation of the current framework to more realistic sea ice simulations that captures this physics is subject to future work.}

There are several ongoing work directions based on this work. One is the coupling with atmosphere models where the atmospheric winds play a role in the sea ice floe dynamics. This is more challenging as the coupling is with the ocean, atmosphere, and sea ice, where the simulation of the dynamics of each is complex. Another future work is data assimilation under this framework. The large-scale model may be used to assimilate the large-scale ocean drag forces, while the local fine-scale DEM may capture the small scales of the ocean eddies. For example, in an ideal setting with a coarse grid, each cell has an eddy localized in the center of the grid box. Then a data assimilation technique may be used to predict these eddies in each cell using the fine-scale DEMs. This would significantly reduce the computational cost, especially in a data assimilation framework that utilizes ensemble Kalman filters.

\bibliographystyle{plain}
%\bibliographystyle{siam}
%\bibliography{ref.bib}
%\bibliographystyle{siamplain}
%\bibliographystyle{siam}
\bibliography{ref}

\end{document}